\documentclass[12pt]{article}
\usepackage{epsfig,wrapfig,fancybox,psboxit,graphics}
\usepackage{amssymb}

\newtheorem{theorem}{Theorem}[section]
\newtheorem{lemma}[theorem]{Lemma}

\newtheorem{definition}[theorem]{Definition}

\newtheorem{cor}[theorem]{Corollary}

\DeclareSymbolFont{AMSb}{U}{msb}{m}{n}
\DeclareMathSymbol{\bdC}{\mathbin}{AMSb}{'103}
\DeclareMathSymbol{\bdR}{\mathbin}{AMSb}{'122}
\newcommand{\bdu}{\mathbf{u}}
\newcommand{\bdv}{\mathbf{v}}
\newcommand{\bdw}{\mathbf{w}}
\newcommand{\bdx}{\mathbf{x}}
\newcommand{\bdy}{\mathbf{y}}
\newcommand{\bdz}{\mathbf{z}}
\newcommand{\bda}{\mathbf{a}}
\newcommand{\bdb}{\mathbf{b}}

\newcommand{\bde}{\mathbf{e}}
\newcommand{\bdf}{\mathbf{f}}
\newcommand{\bdg}{\mathbf{g}}
\newcommand{\bdh}{\mathbf{h}}

\newcommand{\bdr}{\mathbf{r}}
\newcommand{\bdp}{\mathbf{p}}
\newcommand{\bdq}{\mathbf{q}}
\newcommand{\bds}{\mathbf{s}}

\newcommand{\al}{\alpha}

\newcommand{\dl}{\delta}
\newcommand{\Dl}{\Delta}
\newcommand{\eps}{\varepsilon}
\newcommand{\sg}{\sigma}
\newcommand{\vs}{\varsigma}

\newcommand{\Frac}[2]{\mbox{$\displaystyle \frac{#1}{#2}$}}
\newcommand{\pd}[2]{\frac{\partial #1}{\partial #2}}
\newcommand{\bnorm}[1]{\Big\|\, #1 \,\Big\|_2}
\newcommand{\sW}{\mbox{\tiny $W$}}
\newcommand{\qed}{${~} $ \hfill \raisebox{-0.3ex}{\LARGE $\Box$}}

\begin{document}

\title{Computing multiple roots of inexact polynomials}

\author{Zhonggang Zeng \thanks{Department of Mathematics, 
Northeastern Illinois University, Chicago, IL 60625,
email: zzeng@neiu.edu}}

%    General info
%\subjclass{Primary 12Y05, 65H05; Secondary 65F20, 65F22, 65F35}

\date{January 20, 2004.}

%\keywords{polynomial, root, multiplicity}
\maketitle

\begin{abstract}
We present a combination of two algorithms that accurately
calculate multiple roots of general polynomials.
Algorithm I transforms the singular root-finding into a regular nonlinear
least squares problem on a pejorative manifold, and calculates multiple
roots simultaneously from a given multiplicity structure and initial
root approximations.
To fulfill the input requirement of Algorithm I, we develop a numerical
GCD-finder containing a successive singular value updating and an
iterative GCD refinement as the main engine of Algorithm II that 
calculates the multiplicity structure and the initial root approximation.
%
%{\bf %%%%
While limitations of our algorithm exist in identifying the multiplicity
structure in certain situations,
%} %%%%
the combined method calculates multiple  
roots with high accuracy and consistency in practice
without using multiprecision arithmetic even if the
coefficients are inexact.
This is perhaps the first blackbox-type root-finder with such capabilities.
%
%The total complexity of our algorithm is no more than $O(n^3)$ for a
%polynomial of degree $n$.
%
To measure the sensitivity of the multiple roots, a structure-preserving
condition number is proposed and error bounds are established.
%
%{\bf %%%%
According to our computational experiments and error analysis,
%} %%%%
a polynomial being ill-conditioned in the conventional sense can be well
conditioned with the multiplicity structure being preserved,
and its multiple roots can be computed with
high accuracy.
\end{abstract}

\section{Introduction}

In this paper, we present a combination of two numerical 
algorithms for computing multiple roots and the 
multiplicity structures of polynomials.
%
%{\bf %%%%
According to our extensive computational experiments 
%} %%%%
and error estimates, 
the method accurately calculates polynomial roots of non-trivial 
multiplicities without using multiprecision arithmetic, even if 
the coefficients are inexact.

Polynomial root-finding is among the classical problems with longest and 
richest history. 
One of the most difficult issues in root-finding is computing multiple 
roots. 
In addition to requiring {\em exact} coefficients, multiprecision 
arithmetic may be needed when multiple roots are 
present \cite[p. 196]{victorpan97}. 
In fact, using multiprecision has been a common
practice in designing root-finding algorithms and softwares, such 
as those in \cite{bini-mpsolve,far-lou-77,fortune}.
Moreover,  there is a so-called ``attainable accuracy'' in 
computing multiple roots \cite{igarashi-ypma,victorpan97,ypma}: 
to calculate an $m$-fold root to the precision of $k$ correct digits, 
the accuracy of the polynomial coefficients {\em and} the machine precision 
must be at least $mk$ digits. 
This ``attainable accuracy'' barrier also suggests the need of using
multiprecision arithmetic.
Multiprecision softwares such as \cite{bailey} are available. 
However, when polynomial coefficients are truncated, multiple roots 
turn into clusters, and extending machine precision will
never reverse clusters back to multiple roots.
In the absence of accurate methods that are independent of multiprecision 
technology, multiple roots of perturbed polynomials would indeed be 
intractable.

While multiple roots are considered hypersensitive in numerical 
computation, W. Kahan \cite{kahan72} proved 
that if the multiplicities are preserved, those roots may actually 
be well behaved.
More precisely, polynomials with a fixed multiplicity structure form a 
pejorative manifold.
A polynomial is ill-conditioned if it is near such a manifold.
On the other hand, for the polynomial on the pejorative manifold,
its multiple roots are insensitive to multiplicity-preserving perturbations, 
unless the polynomial is also near a submanifold of 
higher multiplicities.
Therefore, to calculate multiple roots accurately, it is important to
maintain the computation on a proper pejorative manifold.

In light of Kahan's theoretical insight, 
we propose Algorithm I in \S \ref{meth2}
that transforms the singular root-finding into a regular nonlinear least 
squares problem on a pejorative manifold. 
By projecting the given polynomial onto the manifold, the computation 
remains structure-preserving. 
As a result, the roots can be calculated simultaneously 
and accurately. 
 
%{\bf %%%%
In applying Algorithm I, one needs a priori knowledge on the 
initial root approximation as well as
multiplicity identification, which is often attempted by
estimation (e.g. \cite{stolan,ypma}) or
clustering (e.g. \cite{brugn,far-lou-85,miyakoda}) 
with unknown certainty.
To fulfill the input requirement of Algorithm I, 
it is preferable to have an algorithm that systematically 
calculates the multiplicity structure.
One of the main difficulties in identifying the multiplicity structure
is the lack of robust method computing the polynomial 
greatest common divisor (GCD).

Recently, many different approaches and strategies have been proposed 
for the numerical
computation of approximate GCD of 
univariate polynomials 
\cite{chin-corless,corless-gianni,emiris-galligo-lombardi,
hribernig-stetter,karmarkar-lakshman,pan96,rupprecht}.
Inspired by those endeavors, especially the pioneering work of 
Corless, Gianni, Trager and Watt
\cite{corless-gianni}, which identifies the GCD degree
by the total SVD (singular value decomposition) of the full 
Sylvester matrix followed by the suggestion of four possible alternative 
ways to compute the GCD using the degree, we propose a numerical 
GCD-finder that employs a successive updating on a sequence of Sylvester
sub-matrices for the smallest singular values only,
followed by
extracting the degree {\em and} the coefficients of the GCD decomposition 
from the singular vector as the initial iterate, and finally applies the
Gauss-Newton iteration to refine the approximate GCD decomposition.
As a result, the GCD-finder is a blackbox-type algorithm in its own right
and constitutes the 
main engine of our proposed Algorithm II in \S \ref{part:gcd}
which, with some limitations specified in \S \ref{lim}, 
calculates the multiplicity structure and its
initial root approximation for a given polynomial.
%} %%%%

In \S \ref{sec:cond}, we propose a structure-preserving condition number 
that measures the sensitivity of multiple roots. 
A polynomial that is ill-conditioned in conventional sense can be well
conditioned with the multiplicity structure being preserved, 
and its roots can be calculated far beyond the barrier of 
``attainable accuracy''. 
This condition number can easily be calculated. 
Error bounds on the roots are established for inexact polynomials. 

In \S \ref{sec:nrII} and \S \ref{sec:nrI}, we present separate numerical 
results for Algorithm I and Algorithm II. 
The numerical results for the combined algorithm are shown in \S \ref{sec:nr}. 
Both algorithms and their combination are implemented as a Matlab package 
{\sc MultRoot} which is electronically 
available from the author. 
This paper also elaborates on certain preliminary results in 
\cite{zeng_mult} by providing detailed proofs and discussions.
 
%{\bf %%%%
Although our emphasis at this stage is on the accuracy
rather than the fast computation, we prove local convergence of our 
Algorithm I when the multiplicity structure and close initial 
approximations are available. We also prove the local convergence
of the GCD iteration in the midst of Algorithm II which
calculates the multiplicity structure with high 
accuracy and consistency in practice along with initial root approximations .  
%} %%%%
%
The combined algorithm 
takes the coefficient vector as the {\em only} input and the outputs results
that include the roots and their multiplicities as well as
the backward error, the estimated forward error, and the structure-preserving 
condition number. 
%
%The total complexity is clearly no more than $O(n^3)$. 
%
The most significant features of the algorithm are its high accuracy 
and its robustness in handling inexact data. 
As shown in numerical examples, the code accurately identifies the 
multiplicity structure and multiple roots for polynomials with a coefficient
accuracy being as low as 7 digits. 
With given multiplicities, Algorithm I converges  even with data 
accuracy as low as 3 decimal digits.
%
%{\bf %%%%
While limitations exist when the polynomial is ill-conditioned
in the sense of structure-preserving sensitivity we define
in \S \ref{sec:cond}, 
%} %%%%
the code appears to be the first blackbox-type root-finder with the
capability to calculate roots and multiplicities 
beyond the barrier of ``attainable accuracy''.

While numerical experiments reported in the literature seldom reach
multiplicity ten, we successfully tested our algorithms on polynomials with 
root multiplicities as high as 400 without using multiprecision arithmetic.
We are aware of no other reliable methods that calculate multiple roots 
accurately by using standard machine precision.
Accurate results for multiple root computation we have seen in the 
literature, such as the methods in \cite{far-lou-77}, 
can be repeated only if multiprecision is used on exact 
polynomials.
A zero-finder for general analytic functions with multiple zeros
has been developed by Kravanja and Van Barel \cite{kra-van}.
The method uses an accuracy refinement with modified Newton's iteration
that may also require multiprecision for multiple roots unless the 
polynomial is already factored \cite{zeng-98}.

There exist general-purpose root-finders using $O(n^2)$ flops or less, 
such as those surveyed in \cite{victorpan97}.
However, the barrier of ``attainable accuracy'' may prevent 
those root-finders from calculating multiple roots accurately 
when the polynomials are inexact (e.g. see Figure 10 in \S 4.6) even if
multiprecision is used.
Our algorithms provide an option of reaching high accuracy on
multiple roots at higher computing cost of $O(n^3)$ which 
may not be a lofty price to pay.

The idea of exploiting the pejorative manifold and the problem structure
has been applied extensively for ill-conditioned problems.
Besides Kahan's pioneering work 30 years ago, theories and computational
strategies for the matrix canonical forms have been studied, such as
\cite{demmel-kagstrom,edelman-elmroth-kagstrom_1,edelman-elmroth-kagstrom_2,lippert-edelman},
to take advantage of the pejorative manifolds or varieties.
At present, it is not clear if those methods can be applied to polynomials
with multiple roots. 

\section{Preliminaries} \label{prelim}

\subsection{Notations}
In this paper, $\bdR^n$ and $\bdC^n$ denote the $n$ dimensional real and
complex vector spaces respectively.
Vectors, always considered columns,  are denoted by boldface lower 
case letters and matrices are denoted by upper case letters.
Blank entries in a matrix are filled with zeros.
The notation $(\cdot)^\top$ represents the transpose of $(\cdot)$, and
$(\cdot)^H$ the Hermitian adjoint (i.e. conjugate transpose) of $(\cdot)$.
When we use a (lower case) letter, say $p$, to denote a polynomial of 
degree $n$, then ~$p_0,p_1,\cdots,p_n$~ are its coefficients as in
\[ p(x) = p_0x^n + p_1 x^{n-1} + \cdots + p_n. \]
The same letter in boldface (e.g. $\bdp$) denotes the coefficient (column)
vector
\[ \bdp = (\,p_0, p_1, \cdots, p_n\,)^\top \]
unless it is defined otherwise. 
The degree of $p$ is $deg(p)$. For a pair of polynomials
$p$ and $q$, their greatest common divisor (GCD) is denoted by $GCD(p,q)$.

\subsection{Basic definitions and lemmas}
\label{defslems}

\begin{definition}
Let \
\( p(x) = p_0 x^n + p_1 x^{n-1} + \cdots + p_n \)~
be a polynomial of degree $n$. For any integer $k\geq 0$, the matrix
\[ C_k(p) =
\begin{array}{c}
\overbrace{\mbox{\ \ \ \ \ \ \ \ \ \ \ \ \ \ \ \ \ \ }}^{k+1}
\\
\left[
\begin{array}{ccc}
p_0 & &  \\ p_1 & \ddots & \\
\vdots & \ddots & p_0 \\ p_n & & p_1 \\
& \ddots & \vdots \\ & & p_n
\end{array} \right]
\end{array}
\]
is called the $k$-th order {\bf convolution matrix} associated with $p$.
\end{definition}

\begin{lemma} \label{lem:conv}
~Let $f$ and $g$ be polynomials of degrees $n$ and $m$ respectively with
$h(x) = f(x)g(x)$.
Then ~$\bdh$~ is the {\bf convolution} of ~$\bdf$~ and ~$\bdg$~ defined by
\newline
\centerline{~\( \bdh = conv(\bdf,\bdg) = C_{m}(f) \bdg = C_{n}(g) \bdf. \)}
\end{lemma}

{\em Proof.} A straightforward verification. \qed

\begin{definition} \label{def:syl}
Let $p$ be a polynomial of degree $n$ and $p'$ be its derivative.
For $k = 1, 2, \cdots n-1$, the matrix of size $(n+k) \times (2k+1)$
\[ S_k(p) = \left[ \begin{array}{ccc} C_{k}\left(p'\right)&
\Big|  & C_{k-1}(p) \end{array} \right]
\]
is called the $k$-th {\bf Sylvester discriminant matrix}.
\end{definition}

\begin{lemma} \label{lem:gcdsv}
Let $p$ be a polynomial of degree $n$ and $p'$ be its derivative
with ~$u = GCD(p,p')$.  
For $j=1,\cdots, n$,
let $\vs_j$ be the smallest singular value of $S_j(p)$.  
Then the following are equivalent
\begin{itemize}
\item[\em (a)] \ \ $deg(u) = m$,
\item[\em (b)] \ \ $p$~ has ~$k=n-m$~ distinct roots,
\item[\em (c)] \ \ \
\( \vs_1,\; \vs_2, \; \cdots, \vs_{k-1} > 0, 
\;\;\;\;\;\; \vs_k = \vs_{k+1} 
= \cdots = \vs_n = 0 \).
\end{itemize}
\end{lemma}

{\em Proof.} 
The equivalence between (a) and (b) is trivial to verify, and the assertion 
that (a) is equivalent to (c) is part of 
Proposition 3.1 in \cite{rupprecht}.
\qed

\begin{lemma} \label{lem:uvw} 
Let $p$ be a polynomial of degree $n$ and $p'$ be its derivative
with $u = GCD(p,p')$ and $deg(u) = m = n-k$. 
Let $v$ and $w$ be polynomials that satisfy
\[ u(x)v(x) = p(x), \;\;\; u(x)w(x) = p'(x). \]
Then
\begin{itemize}
\item[\em (a)] ~$v$ and $w$ are co-prime, namely they have no common factors;
\item[\em (b)] the (column) rank of ~$S_k(p)$ is deficient by one;  
\item[\em (c)] the normalized vector 
$\left[ \begin{array}{r} \bdv \\ -\bdw \end{array} \right]$
is the right singular vector of $S_k(p)$ associated with
the smallest (zero) singular value $\vs_k$;
\item[\em (d)] if $\bdv$ is known, the coefficient vector $\bdu$ of 
$u = GCD(p,p')$ is the solution to the linear system
~~\( C_{m}(v) \bdu = \bdp. \)
\end{itemize}
\end{lemma}

{\em Proof.}
Assertion (a) is trivial. 
~\( S_k(p) \left[ \begin{array}{r} \bdv \\ -\bdw \end{array} \right]
= C_{k}\left(p'\,\right)\bdv - C_{k-1}(p) \bdw = 0\)~
because it is the coefficient vector of
~\( p' v - p w \; \equiv \; (uw)v - (uv)w \; \equiv \; 0  \).
Let ~$\hat{\bdv} \in \bdC^{k+1}$ and $\hat{\bdw} \in \bdC^k$~ be 
coefficient vectors of polynomials $\hat{v}$ and $\hat{w}$ respectively 
that also satisfy
~$C_{k}\left(p'\,\right)\hat{\bdv} - C_{k-1}(p) \hat{\bdw} = 0$.
Then we also have ~$(uw)\hat{v} - (uv)\hat{w} = 0$, namely 
~$w\hat{v} = v\hat{w}$. Since ~$v$~ and ~$w$~ are co-prime, 
there is polynomial ~$c$~ such that ~$\hat{v}=cv$~ and ~$\hat{w}=cw$ and $c$
is obviously a constant.
Therefore, the single vector 
~$\left[ \begin{array}{r} \bdv \\ -\bdw \end{array} \right]$~
forms the basis for the null space of $S_k(p)$.
Consequently, both assertions (b) and (c) follow. 
The assertion (d) is a direct consequence of Lemma \ref{lem:conv}.
\qed

\begin{lemma} \label{lem:invit}
~Let ~$A \in \bdC^{n\times m}$~ with $n\geq m$ be a matrix whose smallest 
two distinct singular values are $\hat{\sg} > \tilde{\sg}$.
Let ~$Q\left( \begin{array}{c} R \\ 0 \end{array} \right) = A$~ 
be the QR decomposition of $A$
where $Q \in \bdC^{n\times n}$ is unitary and $R \in \bdC^{m\times m}$ is 
upper triangular.
From any vector ~$\bdx_0 \in \bdC^m$~ that is not orthogonal to the right
singular subspace of $A$ associated with $\tilde{\sg}$,
we generate the sequences $\{\sg_j\}$ and $\{\bdx_j\}$, by the 
inverse iteration
\begin{eqnarray} \label{invit}
\left\{ \begin{array}{lrcll}
\mbox{Solve } & R^H \bdy_j & = & \bdx_{j-1} & \mbox{for \ } \bdy_j \in \bdC^m \\
\mbox{Solve } & R \,\bdz_j & =  &\bdy_j & \mbox{for \ }
\bdz_j \in \bdC^m\\
&& \\
\mbox{Calculate }& \bdx_j &= &\Frac{\bdz_j}{\|\bdz_j\|_2}, &
\sg_j =  %\Frac{\bdy_i}{\|\bdz_j\|_2}  =
\bnorm{R \bdx_j} \end{array} \right. &&   j = 1, 2, \cdots. 
\end{eqnarray} 
Then 
~\(\displaystyle \lim_{j\rightarrow \infty} \sg_j = \lim_{j\rightarrow \infty} 
\|A\bdx_j\|_2 = \tilde{\sg} \)~ and
\[
\sg_j = \bnorm{A\bdx_j } = \tilde{\sg} + O\left( \tau^{j} \right), 
\;\;\;\; \mbox{where \ } \tau = \left(
 \frac{\tilde{\sg}}{\hat{\sg}} \right)^2.
\]
If $\tilde{\sg}$ is simple, then $\bdx_j$ converges to the right singular
vector $\tilde{\bdx}$ of $A$ associated with $\tilde{\sg}$.
\end{lemma}

{\em Proof.}
See \cite{vanhuffel} for straightforward verifications.
\qed

\subsection{The Gauss-Newton iteration}

The Gauss-Newton iteration is an effective method for solving nonlinear
least squares problems.
Let $G\,:\,\bdC^m \longrightarrow \bdC^n$ with $n>m$, and $\bda \in \bdC^n$.
The nonlinear system
$G (\bdz) = \bda$ for $\bdz \in \bdC^m$
is overdetermined with no conventional solutions in general.
We thereby seek a {\em weighted least squares solution}. 
Let $W = diag(\omega_1,\cdots,\omega_n)$ be a diagonal weight matrix with
positive weights $\omega_j$'s. 
Let $\|\cdot\|_W$ denote the weighted 2-norm:
\begin{equation} \label{wnorm}
 \Big\|\,\bdv\,\Big\|_W \equiv \Big\| \,W\,\bdv\,\Big\|_2 \equiv
\sqrt{\sum_{j=1}^n \omega_j^2v_j^2}, \mbox{\ \ for all \ \ }
\bdv=(v_1,\cdots,v_n)^\top \in \bdC^n.
\end{equation}
Our objective is to solve the minimization problem
~\( \displaystyle
\min_{\bdz\in \bdC^m} \Big\| G(\bdz) - \bda \Big\|_W^2  \).

\begin{lemma} \label{lem:perp} \ \
~Let ~$F\;:\;\bdC^m \longrightarrow \bdC^n$ 
be analytic with Jacobian being $\mathcal{J}(\bdz)$.
%whose components are analytic in
%every variable entry of $\bdz$.
%
%Let $\mathcal{J}(\bdz)$ be the Jacobian of $F(\bdz)$.
%
If there is a neighborhood $\Omega$ of ~$\tilde{\bdz}$ in $\bdC^m$ such that
~\(  \|\,F(\tilde{\bdz}) \,\|_2  \leq
\|\, F(\bdz) \,\|_2 \) for all $\bdz \in \Omega$,
then
~\( \mathcal{J}(\tilde{\bdz})^H F(\tilde{\bdz}) = 0. \)
\end{lemma}

{\em Proof.} The real case ~$F\,:\, \bdR^m \longrightarrow \bdR^n$ of the lemma
is proved in \cite{dennis-schnabel}. The proof for the complex case is
nearly identical, except for using the Cauchy-Riemann equation.
\qed

By Lemma \ref{lem:perp}, let $J(\bdz)$ be the Jacobian of $G(\bdz)$.
To find a local minimum of $\|\, F(\bdz) \|_2 \equiv
\| \, W \,[ G(\bdz) - \bda \,] \|_2$
with $\mathcal{J}(\bdz) = WJ(\bdz)$, we look for $\tilde{\bdz}\in
\bdC^m$ satisfying
\[
\mathcal{J}(\tilde{\bdz})^H F(\tilde{\bdz}) =
\Big[ WJ(\tilde{\bdz}) \Big]^H\, W \Big[G (\tilde{\bdz}) - \bda \Big] =
J(\tilde{\bdz})^H\, W^2 \Big[\,G (\tilde{\bdz}) - \bda \, \Big] = 0.
\]
In other words, ~$G(\tilde{\bdz}) - \bda$~ is orthogonal,
with respect to $ \langle \bdv, \bdw \rangle \equiv \bdv^H W^2 \bdw$,
to the tangent plane of the manifold
~$ \Pi = \left\{\, \bdu = G(\bdz) \;\Big|\; \bdz \in \bdC^m
\,\right\} $~ at ~$\tilde{\bdu} = G (\tilde{\bdz})$.

\begin{figure}[ht]
\centerline{ \epsfig{figure=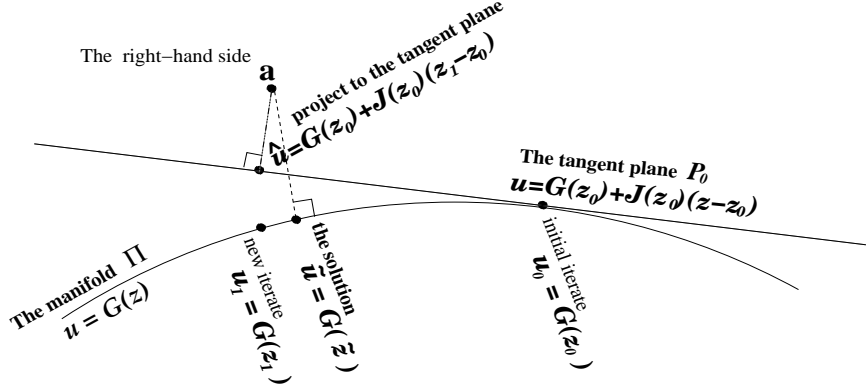,
height=2.1in,width=4.5in} }
\caption[The Gauss-Newton iteration]{Illustration of the 
Gauss-Newton iteration}
\label{odnewton}
\end{figure}

The Gauss-Newton iteration can be derived as follows 
(see Figure \ref{odnewton}).
To find a least squares solution $\bdz = \tilde{\bdz}$ to the 
equation $G(\bdz) = \bda$, 
we look for the point 
$\tilde{\bdu} = G(\tilde{\bdz})$ that is the orthogonal projection of $\bda$
onto $\Pi$.
Let $\bdu_0 = G \left(\bdz_{0}\right)$ in $\Pi$ be near
$\tilde{\bdu} = G(\tilde{\bdz})$. 
We can approximate the manifold $\Pi$ with the tangent plane
\( P_0 = \Big\{\,G \left(\bdz_{0}\right) +
J(\bdz_0)\,(\bdz - \bdz_0)\;\Big|\; 
\bdz \in \bdC^m \; \Big\}  \).
Then the point $\bda$ is orthogonally projected onto the tangent plane $P_0$
at $\hat{\bdu} = G(\bdz_0) + J(\bdz_0)(\bdz_1-\bdz_0)$ by solving the
overdetermined linear system
\begin{equation} \label{wlsqs} G \left(\bdz_{0}\right) +
J(\bdz_0)\,(\bdz - \bdz_0) = \bda \mbox{ \ \ \ or \ \ \ }
J(\bdz_0)\,(\bdz - \bdz_0) = -[G \left(\bdz_{0}\right) - \bda]
\end{equation}
for its weighted least squares solution
\begin{equation} \label{z1}
 \bdz_1 = \bdz_0 - \Big[J(\bdz_0)_{\sW}^+\Big]\,[G(\bdz_0)
- \bda] \mbox{\ \ with \ \ }
J(\bdz_0)_{\sW}^+ =
\Big[ J(\bdz_0)^H W^2J(\bdz_0) \Big]^{-1} J(\bdz_0)^HW^2.
\end{equation}
%
%Here, ~$J(\bdz_0)_{\sW}^+$~ is the weighted pseudo-inverse
%\[ J(\bdz_0)_{\sW}^+ =
%\Big[ J(\bdz_0)^H W^2J(\bdz_0) \Big]^{-1} J(\bdz_0)^HW^2 \]
As long as $J(\bdz_0)$ is of full (column) rank, 
the pseudo inverse ~$J(\bdz_0)_{\sW}^+$~ exists. Therefore
$\bdu_1 = G(\bdz_1)$ is well defined and is expected to be a better 
approximation to
$\tilde{\bdu} = G(\tilde{\bdz})$ than $\bdu_0 = G(\bdz_0)$. 
The Gauss-Newton iteration is then a recursive application of (\ref{z1}) 
(also see \cite{dedieu-shub,dennis-schnabel}).

The convergence theory of the Gauss-Newton iteration has been well
established for overdetermined systems in real spaces
\cite{dennis-schnabel}.
The following lemma is a straightforward generalization of Theorem 10.2.1 in 
\cite{dennis-schnabel} to complex spaces.
Since the lemma itself as well as the proof are nearly identical to those in 
the real case in \cite{dennis-schnabel}, we shall present the lemma
without proof. 

\begin{lemma} \ \ \  \label{conv_lem}
~Let $\Omega \subset \bdC^m$ be a bounded open convex set and
$F\,:\,D \subset \bdC^m \longrightarrow \bdC^n$ be
analytic in an open set $D \supset \overline{\Omega}$.
Let $\mathcal{J}(\bdz)$ be the Jacobian of $F(\bdz)$.
Suppose that there exists $\tilde{\bdz} \in \Omega$ such that
$\mathcal{J}(\tilde{\bdz})^H F(\tilde{\bdz})=0$ with $\mathcal{J}(\tilde{\bdz})$
full rank.
Let $\sg$ be the smallest singular value of $\mathcal{J}(\tilde{\bdz})$.
Let $\dl \geq 0$ be a constant such that
\begin{equation} \label{1024}
\left\|\,\Big[\,\mathcal{J}(\bdz)-\mathcal{J}(\tilde{\bdz})\,\Big]^H F(\tilde{\bdz})\,
\right\|_2 \leq
\dl \, \Big\|\,\bdz-\tilde{\bdz}\,\Big\|_2 \;\;\; \mbox{\ for all \ } 
\bdz \in \Omega.
\end{equation}
If ~$\dl < \sg^2$, then for any ~$c \in \left(\frac{1}{\sg},
\frac{\sg}{\dl}\right)$, there exists $\eps > 0$ such that for all
$z_0 \in \Omega$ with $\| \bdz_0 - \tilde{\bdz} \|_2 < \eps$,
the sequence generated by the Gauss-Newton iteration
\[ \bdz_{k+1} = \bdz_k -  \mathcal{J}(\bdz_k)^+
F(\bdz_k),
\mbox{\ \ where \ } 
\mathcal{J}(\bdz_k)^+ = [\mathcal{J}(\bdz_k)^H\mathcal{J}(\bdz_k)]^{-1}
\mathcal{J}(\bdz_k)^H
\]
for $k=0,1,\cdots$~ is well defined inside 
$\Omega$, converges to $\tilde{\bdz}$, and satisfies
\begin{equation} \label{1025}
 \Big\| \bdz_{k+1} - \tilde{\bdz} \Big\|_2 \leq \frac{c\dl}{\sg} \,
\Big\| \bdz_k - \tilde{\bdz} \Big\|_2 + \frac{c\al \gamma}{2\sg} \,
\Big\| \bdz_k - \tilde{\bdz} \Big\|_2^2,
\end{equation}
%\begin{equation} \label{1026}
% \Big\| \bdz_{k+1} - \tilde{\bdz} \Big\|_2 \leq \frac{c\dl+\sg}{2\sg} \, \Big\|
%\bdz_{k} - \tilde{\bdz} \Big\|_2 < \Big\| \bdz_{k} - \tilde{\bdz} \Big\|_2,
%\end{equation}
where $\al>0$ is the upper bound of
~$\|\mathcal{J}(\bdz)\|_2 $~ on ~$\overline{\Omega}$,
and $\gamma>0$ is the Lipschitz constant of $\mathcal{J}(\bdz)$ in $\Omega$,
namely,
~\( \|\,\mathcal{J}(\bdz+\bdh)-\mathcal{J}(\bdz) \, \|_2  \leq \gamma \, 
\| \,  \bdh \,\| $~
for all $\bdz,\;\bdz+\bdh\;\in\;\Omega . \)
\end{lemma}

\section{Algorithm I: root-finding with given multiplicities}
\label{meth2}

In this section, we assume that the multiplicity structure of a given 
polynomial is known. 
We shall deal with the problem of determining this multiplicity structure
in \S \ref{part:gcd}.
A 
%structure-preserving 
condition number will be introduced to 
measure the sensitivity of multiple roots. 
When the condition number is moderate, the multiple roots can 
be calculated accurately by our algorithm.

\subsection{Remarks on previous work} 

%{\bf %%%%
In Part II of the Technical Report \cite{kahan72}, 
Kahan discussed the sensitivity of
polynomial roots with enlightening insight, and pointed out that it 
may be a misconception to consider multiple roots as infinitely 
ill conditioned. 
Kahan's work on roots of a polynomial $p$ in \cite{kahan72} can be
briefly summarized as follows.
First, Kahan describes the ``pejorative manifold'' of
polynomials with multiple roots. Secondly,
the differentiability is proved for an $m$-fold isolated root with 
respect to coefficients that are constrained to preserve the 
multiplicity $m$ of that root.
This differentiability then naturally leads to the existence of a 
finite local condition number of an {\em isolated} $m$-fold root under
the perturbation that is constrained to preserve the multiplicity $m$ 
of that root. 
Kahan also proves the existence of a vanishing point for 
$p^{(m-1)}(x)$ in a region containing a cluster of $m$ roots of $p$.
Finally, a possible approach is proposed, based on the Lagrange 
multipliers, to find the polynomial nearest to $p$ while possessing 
a $m$-fold root. 

Kahan's work \cite{kahan72} emphasizes on theoretical analysis rather
than computational methodology. The sensitivity analysis is rigorous
while the description of the pejorative manifold is heuristic.
The condition number defined in \cite{kahan72} exists in theory with 
unknown practical attainability.
The implementability of the proposed Lagrange multiplier method in 
numerical computation is still unknown.

In this section, we shall attempt to formulate the pejorative manifold
rigorously. 
More importantly, we emphasize on the practical
computation of multiple roots. 
Our main contribution in this section also include the following:
First, we convert the singular root-finding problem to a least squares
problem and prove its regularity. 
Secondly,
we establish the local convergence of the Gauss-Newton iteration 
solving the least squares problem.
Our third contribution is the formulation of a global
structure-preserving condition number measuring
the combined sensitivity of all roots that are constrained in a
multiplicity structure, instead of an isolated multiple root considered
by Kahan. Not only is this condition number easily computable,
but also it enables us to estimate the computing error
quite accurately according to our error analysis and numerical 
experiment.
Finally, we establish practical procedures that carry out the necessary 
computation on the pejorative manifold. Assembling these elements together, 
we construct our Algorithm I.
%} %%%%%

\subsection{The pejorative manifold} \label{sec:pej}

A polynomial of degree $n$ corresponds to a vector 
(or point) in $\bdC^n$
\[ p(x) =\; p_0x^n + p_1\,x^{n-1} + \cdots + p_n \;\; \sim \;\;
\bda \;=\; (a_1,\cdots,a_n)^\top \equiv \left( \frac{p_1}{p_0},\, \cdots,\,
\frac{p_n}{p_0} \right)^\top,
\]
where ``$\sim$'' denotes this correspondence.
For a partition of $n$, namely a fixed array of positive integers 
$\ell_1, \cdots, \ell_m$ with $\ell_1+ \cdots+ \ell_m=n$, 
a polynomial $p$ that has roots $z_1, \cdots, z_m$ with multiplicities
$\ell_1, \cdots, \ell_m$ respectively can be written as
\begin{equation} \label{expand}
\frac{1}{p_0} p(x) = \prod_{j=1}^m (x-z_j)^{\ell_j} \;\; = \;\; x^n + \sum_{j=1}^n
g_j(z_1,\cdots,z_m)\; x^{n-j},
\end{equation}
where each $g_j$ is a polynomial in $z_1,\cdots,z_m$. 
We have the correspondence
\begin{equation} \label{Gl} p \sim G_\ell (\bdz) \;\equiv\;
\left( \begin{array}{c}
g_1(z_1,\cdots,z_m)\\ \vdots \\ g_n(z_1,\cdots,z_m) \end{array} \right)
\in \bdC^n, \;\;\;
\mbox{\ \ where \ \ }
\bdz = \left( \begin{array}{c} z_1 \\ \vdots \\ z_m \end{array}
\right) \in \bdC^m.
\end{equation}
We now define the pejorative manifold rigorously based on Kahan's
heuristic description.

\begin{definition} \ \
An ordered array of positive integers
$\ell = [ \ell_1,\cdots,\ell_m ] $ is called a {\bf \em
multiplicity structure} of degree $n$ if $\ell_1+\cdots+\ell_m = n$.
For any such given multiplicity structure $\ell$,
the collection of vectors
~\(
 \Pi_\ell \equiv
\left\{\, G_\ell (\bdz) \; \Big| \; \bdz \in \bdC^m \;\right\} \subset \bdC^n
\) 
is called the {\bf \em pejorative manifold} of multiplicity structure
$\ell$, where ~$G_\ell \, : \, \bdC^m \longrightarrow \bdC^n$
defined in (\ref{expand}) -- (\ref{Gl}) is called the {\bf \em
coefficient operator}
associated with the multiplicity structure $\ell$.
\end{definition}

\begin{figure}[htbp]
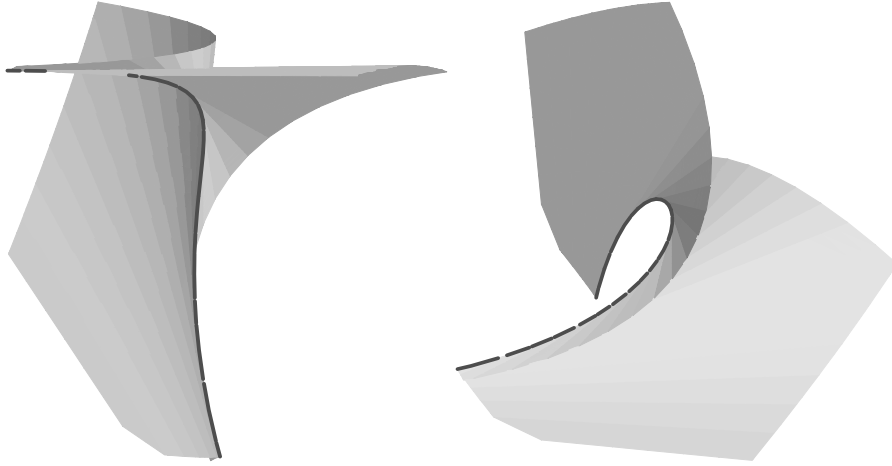

\centerline{
\epsfig{figure=pejo02.eps,height=2.4in,width=2.3in}
\epsfig{figure=pejo03.eps,height=2.4in,width=2.3in}
}
\caption[Pejorative manifolds]{Pejorative
manifolds of polynomials with degree 3 (view from two
angles)}
\label{pej_surf}
\end{figure}

For example, we consider polynomials of degree 3. First, for multiplicity
structure $\ell = [1,2]$, 
\[ (x-z_1)(x-z_2)^2 = x^3 
+ (\mbox{\psboxit{box .8 setgray fill}{$\;-z_1-2z_2\;$}})\,x^2
+ (\mbox{\psboxit{box .8 setgray fill}{$\;2z_1z_2+z_2^2\;$}})\,x 
+ (\mbox{\psboxit{box .8 setgray fill}{$\;-z_1z_2^2\;$}}). \]
A polynomial with one simple root $z_1$ and one double root $z_2$
corresponds to the vector
\begin{equation} \label{3pejcoef}
G_{[ 1,2 ]} (\bdz) \equiv
\left( \begin{array}{r} 
%-z_1-2z_2 
\mbox{\psboxit{box .8 setgray fill}{$\;-z_1-2z_2\;$}}
\\ 
%2z_1z_2+z_2^2 
\mbox{\psboxit{box .8 setgray fill}{$\;2z_1z_2+z_2^2\;$}}
\\
%-z_1z_2^2 
\mbox{\psboxit{box .8 setgray fill}{$\;-z_1z_2^2\;$}}
\end{array} \right) \in \bdC^3, \;\;\;
\mbox{with \ \ }
\bdz  = \left(\begin{array}{c} z_1 \\ z_2 \end{array} \right) \in \bdC^2.
\end{equation}
The vectors  $G_{[1,2]}(\bdz)$ in (\ref{3pejcoef}) for all $\bdz \in 
\bdC^2$ form the pejorative manifold $\Pi_{[1,2] }$.
Similarly, 
\[ \Pi_{[3]} = \Big\{\;
(-3\,z,\,3\,z^2,\,-z^3)^\top \; \Big| \; z \in \bdC \;\Big\}
\]
when $\ell = [3]$. 
$\Pi_{[3]}$ is a submanifold of $\Pi_{[1,2]}$ that contains all polynomials 
with a single triple root.
Figure \ref{pej_surf} shows the manifolds $\Pi_{[1,2] }$
(the wings) and $\Pi_{[3]}$ (the sharp edge) in $\bdR^3$.

As a special case, $\Pi_{[ 1, 1, \cdots,
1]} = \bdC^n$ is the vector space of all monic polynomials
with degree $n$.

\subsection{Solving the nonsingular least squares problem} \label{sec:obj}

Let $\ell = [ \ell_1,\cdots,\ell_m ]$ be a multiplicity structure of
degree $n$ and $\Pi_\ell$ be the corresponding pejorative manifold.
If the polynomial $p \sim \bda \in \Pi_\ell$, then there is a vector
$\bdz \in \bdC^m$ such that $G_\ell(\bdz) = \bda$.
In general, the polynomial system
\begin{equation}
\left\{ \begin{array}{ccc}
       g_1(z_1,\cdots,z_m) & = & a_1 \\
       g_2(z_1,\cdots,z_m) & = & a_2 \\
           \vdots & \vdots & \vdots \\
       g_n(z_1,\cdots,z_m) & = & a_n
     \end{array} \right.
\mbox{ \ \ \ \ or \ \ \ \ }
 G_\ell (\bdz) = \bda %\in \bdC^n, \;\;\;\; \bdz \in \bdC^m
\label{gzal}
\end{equation}
is overdetermined except for the plain structure
$\ell = [ 1, 1, \cdots, 1]$.
Let the weight matrix be $W = diag(\omega_1,\cdots,\omega_n)$.
%with weights $\omega_j>0$, $j=1,\cdots n$, 
and $\|\cdot \|_W$ denote the weighted 2-norm defined in (\ref{wnorm}).
We seek a {\em weighted least squares solution} to (\ref{gzal})
by solving the minimization problem
\begin{eqnarray} 
\lefteqn{
\min_{\bdz\in \bdC^m} \Big\| G_\ell(\bdz) - \bda \Big\|_W^2  \;\;\equiv\;\;
\min_{\bdz\in \bdC^m} \bigg\| W\Big(G_\ell(\bdz) - \bda \Big) \bigg\|_2^2}
\nonumber \\
&& \;\;\equiv\;\; \min_{\bdz\in\bdC^m} \left\{\, \sum_{j=1}^n
\omega_j^2\Big|g_j(\bdz) - a_j\Big|^2 \right\}.\label{wttminl}
\end{eqnarray}
%
%The residual ~$G_\ell(\bdz) - \bda $~ is a by-product of the computation.
%
%Its magnitude is also a measure of the backward error and is handily
%verifiable.

Two common types of weights can be used.
To minimize the overall backward error of the roots, we set
$W = diag(1,1,\cdots,1)$.
On the other hand, the weights
\begin{equation} \label{relweight}
 \omega_j = \min\left\{\,1,\;|a_j|^{-1} \right\}, \;\;\; j = 1, \cdots, n
\end{equation}
lead to minimization of the relative backward error at every coefficient
larger than one.
All our numerical experiments for Algorithm I are conducted using the
weights (\ref{relweight}). 

From Lemma \ref{lem:perp}, let $J(\bdz)$ be the Jacobian of $G_\ell(\bdz)$.
In order to find a local minimum point of $F(\bdz) \equiv W\Big[ 
G_\ell(\bdz) - \bda \Big]$ with $\mathcal{J}(\bdz) = WJ(\bdz)$,
we look for $\tilde{\bdz}\in \bdC^m$ such that 
\begin{equation} \label{perp_eql}
\mathcal{J}(\tilde{\bdz})^H F(\tilde{\bdz}) =
\Big[ WJ(\tilde{\bdz}) \Big]^H\, W \Big[G_\ell (\tilde{\bdz}) - \bda \Big] =
J(\tilde{\bdz})^H\, W^2 \Big[\,G_\ell (\tilde{\bdz}) - \bda \, \Big] = 0.
\end{equation}
%
%Namely, $G_\ell(\tilde{\bdz}) - \bda$~ is orthogonal to the tangent plane
%\[ 
%\left\{\, \bdu = G_\ell(\tilde{\bdz}) + J(\tilde{\bdz}) \,
%(\bdz - \tilde{\bdz}) 
%\in \bdC^n \;\Big|\; \bdz \in \bdC^m \,\right\}
%\]
%of the pejorative manifold
%~$ \Pi_\ell = \left\{\, \bdu = G_\ell(\bdz) \;\Big|\; \bdz \in \bdC^m
%\,\right\}
%$~
%at ~$\tilde{\bdu} = G_\ell (\tilde{\bdz})
%\in \Pi_\ell$.
%The orthogonality is induced from the inner product 
%$ \langle \bdv, \bdw \rangle \equiv \bdv^H W^2 \bdw$.

\begin{definition}
Let $p \sim \bda$ be a polynomial of degree $n$.
For any given multiplicity structure $\ell$ of the same
degree, the vector $\tilde{\bdz}$
satisfying (\ref{perp_eql})
is called a {\bf \em pejorative root vector} or simply
{\bf \em pejorative root} of $p$
corresponding to the multiplicity structure $\ell$ and weight $W$.
\end{definition}

Our algorithms emanate from the following fundamental theorem by 
which one may convert the singular problem of computing multiple roots 
with standard methods to a regular problem by seeking the least squares 
solution of 
(\ref{gzal}).

\begin{theorem}
\label{thm1} \ \
~Let ~$G_\ell\,:\,\bdC^m \longrightarrow \bdC^n$ be the coefficient operator
associated with a multiplicity structure $\ell = [\ell_1,\cdots,\ell_m]$.
Then the
Jacobian $J(\bdz)$ of $G_\ell(\bdz)$ is of full (column) rank
if and only if the entries of ~$\bdz = (z_1,\cdots,z_m)^\top$ are distinct.
\end{theorem}

{\em Proof.} Let $z_1,\cdots,z_m$ be distinct. To prove $J(\bdz)$ is of full 
(column) rank, or the columns of $J(\bdz)$ are linearly independent, 
write 
%$J(\bdz)=\left(\Frac{\partial g_i(\bdz)}{\partial z_j}\right)$ and its 
the $j$-th column of $J(\bdz)$ as
\( J_j = \left( \Frac{\partial g_1(\bdz)}{\partial z_j},
\cdots, \Frac{\partial g_n(\bdz)}{\partial z_j}
\right)^\top. \)
For $j=1,\cdots,m$, let $q_j(x)$, a polynomial in $x$, be defined as follows,
\begin{eqnarray}
q_j(x) & = & \left(\frac{\partial g_1(\bdz)}{\partial z_j}\right)
x^{n-1} + \cdots +
\left(\frac{\partial g_{n-1}(\bdz)}{\partial z_j}\right) x +
\left( \frac{\partial g_n(\bdz)}{\partial z_j} \right) \nonumber \\
& = & \frac{\partial}{\partial z_j} \Big[ x^n + g_1(\bdz)x^{n-1} +
\cdots
+  g_n(\bdz) \Big]
 \; = \;  \frac{\partial}{\partial z_j} \Big[ (x-z_1)^{\ell_1} \cdots
(x-z_m)^{\ell_m} \Big] \nonumber \\
& = & -\ell_j\, (x-z_j)^{\ell_j-1} \left[ \prod_{k\neq j}
(x-z_k)^{\ell_k} \label{jthcol}
\right].
\end{eqnarray}
If \( c_1 J_1 + \cdots + c_m J_m = 0 \)~
for constants $c_1,\cdots,c_m$, then
\begin{eqnarray*}
 q(x) & \equiv & c_1\, q_1(x) + \cdots + c_m\, q_m(x)  
 =  -\sum_{j=1}^m \left\{ c_j \ell_j\, (x-z_j)^{\ell_j-1} \left[
\prod_{k\neq j} (x-z_k)^{\ell_k} \right]\right\} \\
& = &  -\left[\,\prod_{\sg=1}^m (x-z_\sg)^{\ell_\sg-1} \,\right]
\sum_{j=1}^m \left[\,c_j \ell_j
\prod_{k\neq j} (x-z_k) \right] %\right\}
\end{eqnarray*}
is a zero polynomial, yielding
~\(r(x) = \sum_{j=1}^m c_j \ell_j \left[
\prod_{k\neq j} (x-z_k) \right] \equiv 0.
\)
Therefore, for $l = 1,\cdots,m$, 
\( r(z_l) = c_l
\left[\ell_l \prod_{k\neq l} (z_l-z_k) \right] = 0
\)
~implies
~$c_l=0$~ since $\ell_l$'s are positive and $z_k$'s are distinct.
Therefore, $J_j$'s are linearly independent.

On the other hand, suppose $z_1,\cdots,z_m$ are not distinct, say, for 
instance, $z_1=z_2$. 
Then the first two columns of
$J(\bdz)$ are coefficients of polynomials $h_1(x)$ and $h_2(x)$ defined as
\[
-\ell_1\, (x-z_1)^{\ell_1-1}(x-z_2)^{\ell_2}
\prod_{k=3}^m (x-z_k)^{\ell_k}
\mbox{\ \ and \ \ }
-\ell_2\, (x-z_1)^{\ell_1}(x-z_2)^{\ell_2-1} 
\prod_{k=3}^m (x-z_k)^{\ell_k} 
\]
respectively. Since $z_1=z_2$, these two polynomials differ by
constant multiples $\ell_1$ and $\ell_2$. Therefore $J(\bdz)$ is
(column) rank deficient.
\qed

With the system (\ref{gzal}) being nonsingular from Theorem \ref{thm1}, 
the Gauss-Newton iteration
\begin{equation} \label{G-N_it2}
\bdz_{k+1} = \bdz_k - \Big[J(\bdz_k)_{\sW}^+\Big]\,[G_\ell (\bdz_k)-\bda],
\;\;\; k=0,1,\cdots
\end{equation}
on $\Pi_\ell$ is well defined. 
Moreover, 
%When the pejorative root $\tilde{\bdz} =
%(\tilde{z}_1,\cdots,\tilde{z}_m)^\top\in \bdC^m$ has distinct
%components,
%the system (\ref{gzal}) is nonsingular by Theorem \ref{thm1},
%making the Gauss-Newton iteration well defined.
%
we have the convergence theorem
based on Lemma \ref{conv_lem}. 

\begin{theorem} \label{prop2} \ \
~Let ~$\tilde{\bdz} = (\tilde{z}_1,\cdots,\tilde{z}_m)^\top
\in \bdC^m$~ be a
pejorative root of ~$p \sim \bda$~ associated with
multiplicity structure $\ell$ and weight $W$.
Assume $\tilde{z}_1,\, \tilde{z}_2,\,\cdots,\,\tilde{z}_m$ are distinct.
Then there are $\eps,\epsilon >0$ such that, 
if ~$\|\,\bda-G_\ell(\tilde{\bdz})\,\|_W < \eps$~ and
~$\|\, \bdz_0 -\tilde{\bdz}\,\|_2 < \epsilon$, 
the iteration (\ref{G-N_it2}) is well defined and converges to the 
pejorative root $\tilde{\bdz}$ with at least a linear rate.
If we have ~$\bda = G_\ell(\tilde{\bdz})$~ in addition, then the 
convergence is quadratic.
\end{theorem}

{\em Proof.} Let $F(\bdz) = W\Big[G_\ell(\bdz)-\bda\Big]$ and $\mathcal{J}(\bdz)$ 
be its Jacobian.
$F(\bdz)$ is obviously analytic.
From Theorem \ref{thm1}, the smallest singular value $\sg$ of
$\mathcal{J}(\tilde{\bdz})$ is strictly positive.
If $\bda$ is sufficiently close to $G_\ell(\tilde{\bdz})$, then
$\bnorm{F(\tilde{\bdz})} = \Big\|\,G_\ell(\tilde{\bdz})-\bda\,\Big\|_W$~
will be small enough, making (\ref{1024}) holds with $\dl < \sg^2$.
Therefore all conditions of Lemma \ref{conv_lem} are satisfied and there
is a neighborhood $\Omega$
of $\tilde{\bdz}$ such that if $\bdz_0 \in \Omega$,
the iteration (\ref{G-N_it2}) converges and satisfies
(\ref{1025}).
If in addition $\bda = G_\ell(\tilde{\bdz})$, then $F(\tilde{\bdz})=0$
and therefore $\dl=0$ in (\ref{1024}) and (\ref{1025}).
The convergence becomes quadratic. \qed

As a special case for the structure $\ell = [1,1,\cdots,1]$,
equations in (\ref{gzal}) form Vi\'ete's system of $n$-variate polynomial
system.  Solving this system via Newton's iteration is equivalent to
the Weierstrass (Durand-Kerner) algorithm \cite{victorpan97}.
When a polynomial has multiple roots, Vi\'ete's system becomes
singular at the non-distinct root vector. This singularity 
appears to be the very reason that causes the ill-conditioning of 
conventional root-finders: a wrong pejorative 
manifold is used.

\subsection{The structure-preserving condition number}
\label{sec:cond}
There are many insightful discussions on the numerical condition of 
polynomial roots in the literature such as
\cite{demmel-ill,Gautschi-84,kahan72,wilkinson-63,stetter-cond,winkler-01}.
%\cite{demmel-ill,Gautschi-84,kahan72,bezIII,wilkinson-63,stetter-cond,winkler-01}.
%
In general, a condition number can be characterized as the smallest number 
satisfying
\begin{equation} \label{cond}
 \Big[ forward\_error \Big] \leq \Big[ condition\_number \Big] \times
\Big[ backward\_error \Big] + h.o.t., \end{equation} 
where $h.o.t$ represents higher order terms of the backward error.
For a polynomial with multiple roots, under {\em unrestricted} perturbation, 
the only condition number satisfying (\ref{cond}) is infinity. 
For a simple example, let polynomial ~$p(x) = x^2$. 
~A backward error ~$\eps$~ makes the perturbed polynomial 
~$\tilde{p}(x) = x^2+\eps$~ having roots ~$\pm \sqrt{\eps}i$~ with 
forward error ~$\sqrt{\eps}$~ in magnitude. 
The only ``constant'' ~$c$~ which accounts for ~$\sqrt{\eps} \leq c \, \eps$~
for {\em all} ~$\eps>0$~ must be infinity.

By changing the computational objective from solving a polynomial equation
$p(x)=0$ to the nonlinear least squares problem in the form of 
(\ref{wttminl}), the structure-altering noise is filtered out, and 
the multiplicity structure is preserved.
With this shift in computing strategy, 
the sensitivity of the roots can be analyzed differently.

Let's consider the root vector $\bdz$ of 
~$p \sim \bda=G_\ell(\bdz)$. 
The polynomial $p$ is perturbed, 
with multiplicity structure $\ell$ being preserved,
to be ~$\hat{p} \sim \hat{\bda}=G_\ell(\hat{\bdz})$. 
In other words, both $p$ and $\hat{p}$ are on the same pejorative
manifold $\Pi_\ell$.
Then
\[ \hat{\bda}-\bda\ = G_\ell(\hat{\bdz}) - G_\ell(\bdz) = 
J(\bdz)(\hat{\bdz}-\bdz) +
O\left(\|\,\hat{\bdz}-\bdz\,\|^2\right) \]
where $J(\bdz)$ is the Jacobian of $G_\ell (\bdz)$.
Assuming the entries of $\bdz$ are distinct, by Theorem \ref{thm1}.
$J(\bdz)$ is of full rank. 
Consequently,
\begin{eqnarray} 
 \Big\|\,W(\hat{\bda} - \bda)\,\Big\|_2 & = & 
\Big\|\,[WJ(\bdz)](\hat{\bdz}-\bdz)
\Big\|_2+h.o.t.\,, 
\nonumber \\ 
\mbox{namely,\ \ \ \ \ \ } 
\Big\|\,\hat{\bda}-\bda\,\Big\|_W & \geq & \sigma_{min} \Big\|\,
\hat{\bdz}-\bdz\,\Big\|_2 + h.o.t.,  
\nonumber \\
\mbox{or \ \ \ \ \ \ } 
\Big\|\,\hat{\bdz}-\bdz\,\Big\|_2 & \leq & \left(\frac{1}{\sigma_{min}}
\right)\, 
\Big\|\,\hat{\bda}-\bda\, \Big\|_W + h.o.t. \label{est1}
\end{eqnarray}
where ~$\sigma_{min}$~, the smallest singular value of $WJ(\bdz)$, is 
strictly positive since $W$ and $J(\bdz)$ are of full rank.
The distance $\|\,\hat{\bdz}-\bdz\,\|_2$ is the forward error and
the weighted distance $\|\,\hat{\bda}-\bda\, \|_W $ measures the
backward error.
Therefore, the sensitivity of the root vector is asymptotically
bounded by $\frac{1}{\sigma_{min}}$ times the size of the
multiplicity-preserving
perturbation. In this sense,
the multiple roots are not infinitely sensitive.

\begin{definition} \ \ Let $p$ be a polynomial and
$\bdz$ be its pejorative root
corresponding to a given multiplicity structure $\ell$ and weight $W$. 
Let $G_\ell$ be the coefficient operator associated with $\ell$, $J$ 
be its Jacobian, and $\sg_{min}$ be the smallest singular value of 
~$WJ(\bdz)$.
Then the {\bf \em condition number of} ~$\bdz$~ {\bf with respect
to the multiplicity structure} $\ell$ {\bf and weight $W$} is
defined as
\[ \kappa_{\ell,w}(\bdz) = \frac{1}{\sigma_{min}}. \]
\end{definition}

{\bf Remark:} The condition number $\kappa_{\ell,w}(\bdz)$ is 
structure dependent. 
The array $\ell =[\ell_1,\cdots,\ell_m]$ may or may not be the 
{\em actual} multiplicity structure. 
A polynomial has different condition numbers 
corresponding to different pejorative roots on 
various pejorative manifolds.
For example, see Table \ref{table_ex1} in \S \ref{sec:trirt}.

We now estimate the error on pejorative roots of polynomials 
with inexact coefficients.
In this case, the given polynomial $\hat{p}$ is assumed to be
arbitrarily perturbed from $p$ with both polynomials near a
pejorative manifold $\Pi_\ell$. 
In exact sense, neither polynomial possesses the structure $\ell$.
The nearby pejorative manifold causes both polynomials
being ill-conditioned in conventional sense. 
Consequently, the exact roots of $\hat{p}$ can be far from those of 
$p$ even if two polynomials are close to each other. 
However, the following theorem ensures that their pejorative 
roots, not exact roots, may still be insensitive to perturbation. 

\begin{theorem} \label{prop3} 
~For a fixed $\ell = [\ell_1,\cdots,\ell_m]$, 
let the polynomial ~$\hat{p} \sim \hat{\bdb}$~ be
an approximation to ~$p  \sim \bdb$ with pejorative
roots $\hat{\bdz}$ and $\bdz$ respectively that are corresponding
to the multiplicity structure $\ell$ and a weight $W$. 
Assume the components of ~$\bdz$ are distinct while 
$\|\, G_\ell(\hat{\bdz}) - \hat{\bdb} \, \|_W$ reaches a local minimum
at $\hat{\bdz}$. 
If ~$\|\,\bdb - \hat{\bdb} \,\|_W$ and 
~$\|\, G_\ell(\bdz) - \bdb \,\|_W$ are sufficiently small, then
\begin{equation} \label{est3}
\bnorm{\bdz-\hat{\bdz}} \leq 2\cdot \kappa_{\ell,w}(\bdz) 
\cdot %\frac{2}{\sg_{min}}
\left( \Big\|\, G_\ell(\bdz) - \bdb \,\Big\|_W + \Big\|\, \bdb - \hat{\bdb} \,
\Big\|_W \right) + h.o.t.
\end{equation}
\end{theorem}

{\em Proof.} From (\ref{est1}),
\begin{eqnarray*}
\lefteqn{\bnorm{\bdz-\hat{\bdz}}  \leq  \kappa_{\ell,w}(\bdz)\,
\Big\|\, G_\ell(\bdz) - G_\ell(\hat{\bdz}) \, \Big\|_W + h.o.t.}\\
& \leq & \kappa_{\ell,w}(\bdz) \left(\,
\Big\|\, G_\ell(\bdz) -  \bdb \, \Big\|_W +
    \Big\|\, \bdb - \hat{\bdb} \, \Big\|_W +
 \Big\|\, G_\ell(\hat{\bdz}) - \hat{\bdb} \, \Big\|_W  \right)+h.o.t.
\end{eqnarray*}
Since $\|\, G_\ell(\hat{\bdz}) - \hat{\bdb} \, \|_W$ is a local minimum,
we have
%\begin{eqnarray*}
\[ \Big\|\, G_\ell(\hat{\bdz}) - \hat{\bdb} \, \Big\|_W  \leq
\Big\|\, G_\ell(\bdz) - \hat{\bdb} \, \Big\|_W
  \leq
\Big\|\, G_\ell(\bdz) - \bdb \,\Big\|_W + \Big\|\, \bdb - \hat{\bdb} \, 
\Big\|_W,  \]
%\end{eqnarray*}
%
and the assertion of the theorem follows. \qed

By the above theorem, when a polynomial is perturbed, 
the error on the pejorative roots depends on the magnitude of the 
perturbation $\|\, \bdb - \hat{\bdb} \, \|_W$, the distance 
to the pejorative manifold (namely $\|\, G_\ell(\bdz) - \bdb \, \|_W$), 
as well as the condition number $\kappa_{\ell,w}(\bdz)$. 
Although the (exact) roots may be hypersensitive, their
pejorative roots are stable if $\kappa_{\ell,w}(\bdz)$ is moderate.

%If the polynomial $\hat{p}$ is obtained by an arbitrarily perturbation 
%from $p$ that keeps the multiplicity structure $\ell$, 
%the error estimate on the pejorative root of $\hat{p}(x)$ is even simpler.
For a polynomial $p$ having multiplicity structure $\ell$, 
we can now estimate the error of its multiple roots computed from 
its (inexact) approximation $\hat{p}$. The perturbation from $p$ to 
$\hat{p}$ can be arbitrary, such as rounding up digits in coefficients. 
The (exact) roots of $\hat{p}$ are all simple in general and far from
the multiple roots of $p$.  
The following corollary ensures that the {\em pejorative} root 
$\hat{\bdz}$ of 
$\hat{p}$ with respect to the multiplicity structure $\ell$ can 
be an accurate approximation to the multiple roots $\bdz$ of $p$. 

\begin{cor} \label{spccor}
Under the condition of Theorem \ref{prop3}, if
~$\bdz$~ is the exact root vector of $p$ with
multiplicity structure $\ell$, then
\begin{equation} \label{est2}
 \bnorm{\bdz-\hat{\bdz}} \leq 2\cdot \kappa_{\ell,w}(\bdz)\cdot 
%\frac{2}{\sg_{min}}
 \Big\|\, \bdb - \hat{\bdb} \, \Big\|_W  + h.o.t.
\end{equation}
\end{cor} 

{\em Proof.} Since $\bdz$ is exact, $\|\, G_\ell(\bdz) - \bdb \,\|_W = 0$ 
in (\ref{est3}). \qed

The ``attainable accuracy'' barrier suggests that when multiplicity 
increases, the roots sensitivity intensifies. 
However, there is no apparent correlation between the magnitude of the 
multiplicities and the structure-constraint sensitivity.
For example, consider polynomials
\[ p_\ell (x) = (x+1)^{\ell_1} (x-1)^{\ell_2} (x-2)^{\ell_3} \]
\begin{wrapfigure}[9]{l}{2.1in}
\footnotesize \vspace{-6.5mm}
\begin{center}
\begin{tabular}{|rrr|c|} \hline
\multicolumn{3}{|c|}{multiplicities} & condition \\
$\ell_1$ & $\ell_2$ & $\ell_3$ &  number \\ \hline
1 & 1 & 1 & 3.1499 \\
1 & 2 & 3 & 2.0323 \\
10 & 20 & 30 & 0.0733 \\
100 & 200 & 300 & 0.0146  \\ \hline
\end{tabular}
\end{center} \vspace{-3mm}
\caption[The sensitivity and root multiplicities]{
\footnotesize The sensitivity and multiplicities}
\label{tab:pejmult}
\end{wrapfigure}
with different multiplicities $\ell = [\ell_1,\ell_2,\ell_3]$.
For the weight $W$ defined in (\ref{relweight}),
Figure \ref{tab:pejmult} lists the condition numbers for
different multiplicities.  
As seen in this example, 
the magnitude of root error can actually be {\em less} than 
that of the data
error when the condition number is less than one. 
The condition theory described above indicates that multiprecision 
arithmetic may {\em not} be a necessity, and the ``attainable accuracy'' 
barrier appears to be highly questionable.

In \S \ref{sec:nrII} and \S \ref{sec:nr}, more examples will show that our 
iterative algorithm indeed reaches the accuracy permissible by the 
condition number $\kappa_{\ell,w} (\bdz)$, which can be calculated
with negligible cost. 
The QR decomposition of the Jacobian $J(\bdz)$ 
is required by the
iteration (\ref{G-N_it2}), and can be recycled to calculate 
$\kappa_{\ell,w} (\bdz)$. 
The inverse iteration in Lemma \ref{lem:invit} is suitable for
finding the smallest singular value.
%
%Since only the size, rather than the accurate digits, of the 
%smallest singular value is required, a few inverse iteration steps are 
%generally enough. 

\subsection{The numerical procedures} \label{sec:proc}

\begin{figure}[ht]
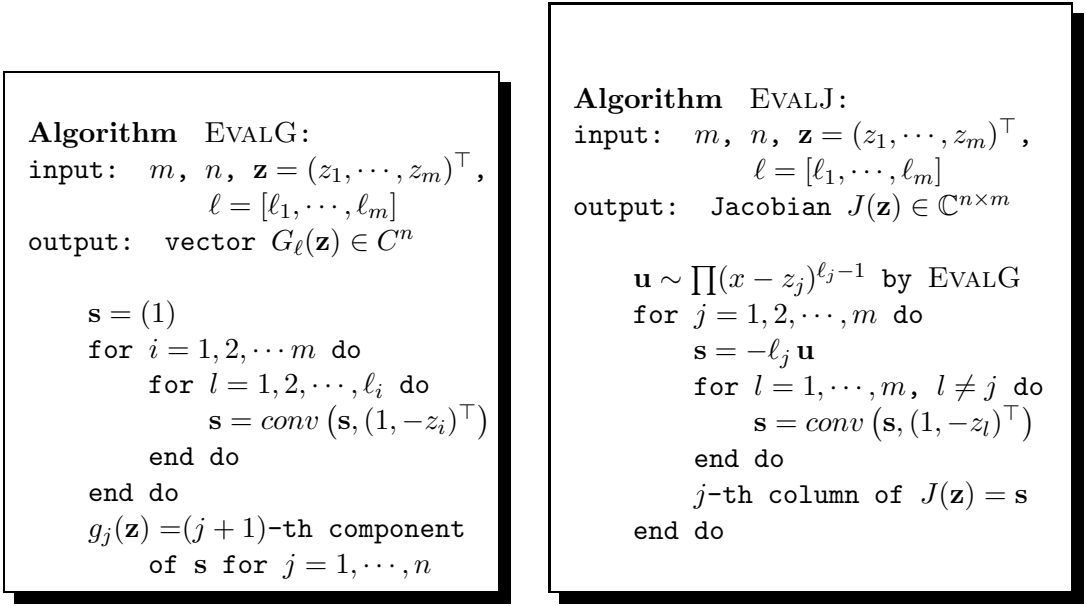

\small  \vspace{-2mm}
\begin{center}
\shadowbox{\parbox{2.0in}{
\baselineskip0mm
{\tt
\begin{tabbing}
\hspace{0mm} \=
\hspace{6mm} \=
\hspace{6mm} \=
\hspace{6mm} \=
\hspace{6mm} \=
\hspace{6mm} \=
\hspace{6mm} \= \\
\>{\bf Algorithm } {\sc EvalG}: \\
\> input: $m$, $n$, $\bdz=(z_1,\cdots,z_m)^\top$, \\
\>\>\>\>$\ell=[\ell_1,\cdots,\ell_m]$ \\
\> output: vector $G_\ell(\bdz) \in C^n$ \\
\>\\
\>\> $\bds = (1) $ \\
\>\> for $i = 1, 2, \cdots m$ do\\
\>\>\> for $l=1,2,\cdots,\ell_i$ do\\
\>\>\>\> $\bds= conv\left(\bds, (1,-z_i)^\top\right)$\\
\>\>\> end do\\
\>\> end do\\
\>\> $g_j(\bdz) = $$(j+1)$-th component \\
\>\>\> of $\bds$ for $j=1,\cdots,n$ \ \ \  %\\
\end{tabbing}
}}} \ \ \ 
\shadowbox{\parbox{5.0in}{
\baselineskip0mm
{\tt
\begin{tabbing}
\hspace{0mm} \= 
\hspace{6mm} \=
\hspace{6mm} \=
\hspace{6mm} \=
\hspace{6mm} \=
\hspace{6mm} \=
\hspace{6mm} \= \\
\>{\bf Algorithm } {\sc EvalJ}: \\
\> input: $m$, $n$, $\bdz=(z_1,\cdots,z_m)^\top$, \ \  \\
\>\>\>\> $\ell=[\ell_1,\cdots,\ell_m]$ \\
\> output: Jacobian $J(\bdz) \in \bdC^{n\times m}$ \\
\>\\ 
\>\> $\bdu \sim \prod (x-z_j)^{\ell_j-1}$ by {\sc EvalG} \\
\>\> for $j=1,2,\cdots,m$ do \\
\>\>\> $\bds = -\ell_j\,\bdu$ \\
\>\>\> for $l = 1, \cdots, m$, $l\neq j$ do\\
\>\>\>\> $\bds= conv\left(\bds,(1,-z_l)^\top\right)$\\
\>\>\> end do\\
\>\>\> $j$-th column of $J(\bdz)=\bds$ \\
%\>\>\>\> components of $\bds$ \\
\>\> end do %\\
%\>
\end{tabbing} 
}}}
\end{center}\vspace{-3mm}
\caption{Pseudo-codes for evaluating $G_\ell(\bdz)$ and
$J(\bdz)$}
\label{alg_gj} \vspace{-2mm}
\end{figure}

The iteration (\ref{G-N_it2}) requires calculation of
the vector value of 
$G_\ell(\bdz_k)$ and matrix value of $J(\bdz_k)$, 
where the components of $G_\ell(\bdz)$ are
defined in (\ref{expand}) and (\ref{Gl}) as coefficients of the polynomial
$p(x) = (x-z_1)^{\ell_1}\cdots (x-z_m)^{\ell_m}.$
While the explicit formulas for each 
$g_j(z_1,\cdots,z_m)$ and $\pd{g_j}{z_i}$
can be symbolically (inefficiently in general) computed using softwares 
like Maple, we propose more 
efficient numerical procedures for computing 
$G_\ell(\bdz)$ and $J(\bdz)$ in Figure \ref{alg_gj}.

The polynomial multiplication
is equivalent to the vector convolution (Lemma \ref{lem:conv}).
The polynomial $p(x) = (x-z_1)^{\ell_1}\cdots (x-z_m)^{\ell_m}$ 
can thereby be constructed 
from recursive convolution with vectors $(1,-z_j)^\top$, $j=1,2,\cdots,m$.
As a result, $G_\ell(\bdz)$ is computed through the nested loops shown in 
Figure \ref{alg_gj} as Algorithm {\sc EvalG}.
It takes $n^2 + O(n)$ flops (additions and
multiplications) to calculate $G_\ell(\bdz)$.

The $j$-th column of the Jacobian $J(\bdz)$,
as shown in the proof of Theorem \ref{thm1},
can be considered the coefficients of the polynomial $q_j(x)$ defined in
(\ref{jthcol}).
The cost of computing $J(\bdz)$ is no more than $mn^2 + O(n)$ flops.
Each step of the Gauss-newton iteration takes $O(nm^2)$ flops. 
Therefore, for a polynomial of degree $n$ with $m$ distinct roots, 
the complexity of Algorithm I is $O(m^2n+mn^2)$. 
The worst case occurs when $m=n$ and the complexity becomes $O(n^3)$.
The complete pseudo-code of the Algorithm I is shown in Figure
\ref{g-n_sc}. 

\begin{figure}[ht]
\small \vspace{-3mm}
\begin{center}
\shadowbox{\parbox{5.7in}{
\baselineskip0mm
\vspace{-2mm}
{\tt
\begin{tabbing}
\hspace{6mm} \=
\hspace{6mm} \=
\hspace{6mm} \=
\hspace{6mm} \=
\hspace{6mm} \=
\hspace{6mm} \=
\hspace{6mm} \= \\ 
\>{\bf Pseudo-code} {\sc PejRoot} ~(Algorithm I): \\
\> ~input: $m$, $n$, $\bda \in \bdC^n$, weight matrix $W$, 
initial iterate $\bdz_0$, \\
\>\>\>\>  multiplicity structure $\ell$, error tolerance $\tau$ \ \  \\
\> ~output: Roots $\bdz=(z_1,\cdots,z_m)$, or message of failure \\
\>\\
\>\> for $k=0,1,\cdots$ do \\
\>\>\> Calculate $G_\ell(\bdz_k)$ and $J(\bdz_k)$ with {\sc EvalG}
and {\sc EvalJ}\\
\>\>\> Compute the least squares solution $\Dl \bdz_k$ to the linear 
\ \ \ \\
\>\>\>\>\> system ~$[WJ(\bdz_k)](\Dl \bdz_k) = W[G_\ell(\bdz_k)-\bda]$ \\
\>\>\> Set $\bdz_{k+1} = \bdz_k - \Dl \bdz_k$ and
$\delta_k = \|\Dl \bdz_k \|_2$ \\
%\> \\
\>\>\> if $k \geq 1$ then \\
\>\>\>\> if $\delta_k \geq \delta_{k-1}$ then, stop, output failure message \\
\>\>\>\> else if $\frac{\delta_k^2}{\delta_{k-1}-\delta_k} <
\tau$ then, \  stop, output $\bdz = \bdz_{k+1}$ \\
\>\>\>\> end if \\
\>\>\> end if \\
\>\> end do %\\
\end{tabbing}
}}}
\end{center} \vspace{-2mm}
\caption[Pseudo-code of Algorithm I]{Pseudo-code of Algorithm I}
\label{g-n_sc} \vspace{-3mm}
\end{figure}

\subsection{Numerical results for Algorithm I} 
\label{sec:nrII}

Algorithm I  is implemented as a Matlab code {\sc PejRoot}.
All the tests of {\sc PejRoot} are conducted with IEEE double precision 
(16 decimal digits) without extension. In comparison, other algorithms and 
software may use unlimited machine precision in some cases. 
%The testing polynomials are in general form.

\subsubsection{The effect of ``attainable accuracy''}

Conventional methods, such as Farmer-Loizou methods \cite{far-lou-77},
are subject to the ``attainable accuracy'' barrier.
We made a straightforward %Matlab and Maple 
implementation of the
Farmer-Loizou third order iteration suggested in \cite{far-lou-77} and
apply it to the same example they used
\[ p_1(x) = (x-1)^4(x-2)^3(x-3)^2(x-4). \]
Both iterations start from $\bdz_0 = (1.1,1.9,3.1,3.9)$
using the standard IEEE double precision.
The ``attainable accuracy'' of the roots are 4, 5, 8, 16 digits
respectively.  For 100 iteration steps,
the Farmer-Loizou method produces iterates that bounce around
the roots.
In contrast, our iteration smoothly converges to the roots and
reaches accuracy of
14 digits. The ``attainable accuracy'' barrier has no effect on our
algorithm. The iterations are shown in Table \ref{tab_farloi1} for three 
roots $x=1,2,3$ with highest multiplicities.

\begin{table}[ht]
\scriptsize
\begin{verbatim}
Farmer-Loizou third order iteration    |                  PejRoot result
step  iterateas                        |step  iterates  
  1  1.0009     1.998      3.001       | 1  1.03              1.8               3.4               
  2  0.99997    1.9999992  3.000000008 | 2  0.997             1.98              2.6               
  3  0.01       3.4        2.9988      | 3  1.00009           2.05              2.8               
  4  0.8        2.3        3.000007    | 4  0.99994           1.994             2.98              
  5  0.998      2.007      3.0000001   | 5  1.000003          2.0001            2.9990            
  6  1.0000007  2.00000007 2.99996     | 6  0.999999997       2.000000005       2.9999990         
...  ...                               | 7  1.00000000000000  2.0000000000002   2.999999999998    
100  1.00000008 3.3        2.99999997  | 8  1.00000000000000  2.00000000000000  2.99999999999999  
\end{verbatim}
\normalsize \vspace{-6mm}
\caption[\small Comparison with the Farmer-Loizou third order iteration in low 
multiplicity case]{
\footnotesize
Comparison with the Farmer-Loizou third order 
iteration in low multiplicity case. Three
roots are shown with unimportant digits truncated}
\label{tab_farloi1}
\vspace{-3mm}
\end{table}

In the same problem, we increase the multiplicities 10 times as large,
generating
\[ p_2(x) = (x-1)^{40}(x-2)^{30}(x-3)^{20}(x-4)^{10} \]
with 16-digit accuracy in coefficients.
In this test, our method still uses the standard
16-digit arithmetic and attains 14 correct digits on the roots, 
while Farmer-Loizou method uses 1000-digit operations in Maple and fails
(Three roots iterations are shown in
Table \ref{tab_farloi2}).

\begin{table}[ht]
\scriptsize
\begin{verbatim}
Farmer-Loizou third order iteration |     PejRoot result            
step  iterateas                     |steps  iterates  
 1   0.47       -.33        3.02    | 1     1.004             1.98              3.05              
 2  32.92      -4.65        2.69    | 2     1.0001            1.998             3.003             
 3   4.75      -1.80        1.75    | 3     0.9999998         2.000006          2.99997           
 4 205.96        .40        1.54    | 4     0.999999999994    2.00000000001     2.9999999990      
...  ...                            | 5     1.00000000000000  2.00000000000001  2.99999999999997  
100  5.99       1.10        0.30    | 6.    1.00000000000000  2.00000000000001  2.99999999999998  
\end{verbatim}
\normalsize \vspace{-4mm}
\caption[Comparison with the Farmer-Loizou third order 
iteration in high multiplicity case]{
\footnotesize
Comparison with the Farmer-Loizou third order 
iteration in high multiplicity case. Three
roots are shown with unimportant digits truncated}
\label{tab_farloi2}
\vspace{-4mm}
\end{table}

The true accuracy barrier of our Algorithm I is the 
condition number $\kappa_{\ell,w}(\bdz)$. 
Matlab constructed the test polynomial with a relative coefficient error of
$4.56\times 10^{-16}$, the condition number is 29.3. 
The root error is approximately $1 \times 10^{-14}$, which is within the 
error bound $2\times(29.3)\times (4.56\times 10^{-16}) = 2.67\times 10^{-14}$ 
established in Corollary \ref{spccor}.

There are state-of-art root-finding packages available
using multiprecision,
such as {\sc MPSolve} implemented by Bini et al \cite{bini-mpsolve} and
{\sc Eigensolve} by Fortune \cite{fortune}. 
%and the Sch\"{o}nhage algorithm implemented by Gourdon \cite{gourdon}. 
%
If the given polynomial is exact (e.g. polynomial with rational coefficients), 
those packages in general are capable of calculating all roots to the desired 
accuracy via extending the machine precision according to the ``attainable
accuracy''. %while multiple roots are still calculated as clusters 
%without multiplicity identification. 
%
For inexact polynomials, the accuracy of those 
packages on multiple roots are limited no matter how many 
digits the machine precision is extended.
For example, consider the polynomial
\[ p(x) = \left(x-\sqrt{2}\right)^{20} \, \left(x-\sqrt{3}\right)^{10}. \]
The coefficients are calculated to 100-digit accuracy. 
The ``attainable accuracy'' for the roots $\sqrt{2}$ and $\sqrt{3}$ are 5 
and 10 digits respectively. {\sc MPSolve} and {\sc Eigensolve} output 
nearly identical results in accordance with this ``attainable accuracy''. 
In contrast, our software using only
16 digits precision in coefficients without extending the machine precision, 
still outputs roots of 15 digits accuracy along with accurate multiplicities 
(see Table \ref{tab_mpsolve}).
%
%Again, the ``attainable accuracy'' has no effect on our algorithm. 
%Moreover, the structure-preserving condition number accurately predicts 
%the root error.

\begin{table}[ht]
\scriptsize
\begin{verbatim}
  MPSolve and Eigensolve results        | MultRoot results with 16-digit input/machine precision
    with 100-digit input accuracy       |
  (unimportant digits are truncated)    | THE CONDITION NUMBER:                     0.90775      
    1.41412    - 0.000013i              | THE BACKWARD ERROR:                 6.66E-016
    1.41412    + 0.000013i              | THE ESTIMATED FORWARD ROOT ERROR:   1.21e-15           
     ... ...                            |
    1.73205077 - 0.0000000094i          | computed roots               multiplicities
    1.73205077 + 0.0000000094i          |      1.732050807568876                   10
     ... ...                            |      1.414213562373096                   20
\end{verbatim}
\normalsize  \vspace{-5mm}
\caption{\footnotesize 
Comparison with multiprecision root-finders {\sc MPSolve} and
{\sc Eigensolve}} 
\label{tab_mpsolve}
\vspace{-6mm}
\end{table}

\subsubsection{Clustered multiple roots}
\label{sec:trirt}

\begin{wrapfigure}[7]{r}{2.5in}
\vspace{-16mm}
\centerline{\epsfig{figure=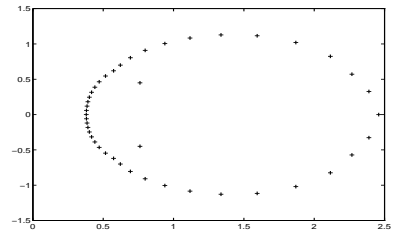,
width=2.00in,height=1.20in}}\vspace{-3mm}
\caption[Roots of a polynomial with 
clustered multiple roots calculated by Matlab]{\footnotesize 
The root cluster from three multiple roots
calculated by Matlab {\tt roots}}
\label{trizero}
\end{wrapfigure}
Let
\newline \( f(x) = (x-0.9)^{18}(x-1)^{10}(x-1.1)^{16}.\)
The roots are highly multiple and clustered.
The Matlab function {\sc roots}
produces $44$ ill-conditioned roots scattered in a box
of $2.0\times 2.0$ (see Figure \ref{trizero}).
In contrast, the Algorithm I code {\sc PejRoot} obtains all three multiple
roots for at least 14 digits in accuracy by taking two additional 
iteration steps on the information of multiplicity structure and the initial
iterate provided by our Algorithm II in \S \ref{part:gcd}.

\centerline{\scriptsize
\vspace{-0mm}
\( \begin{array}{c|clll}
step &&    z_1 & z_2 & z_3 \\ \hline
  0 &&  0.89999999993 & 0.9999999993 & 1.0999999998 \\
  1 &&  0.999999999999991 &      1.00000000000001 & 1.10000000000001   \\
  2 &&  0.999999999999991 &      1.000000000000001 & 1.10000000000001   \\
\end{array}
\)}  
\vspace{2mm}

The backward accuracy can easily be verified to be
less than $1.36 \times 10^{-15}$. 
%Namely, the coefficients of 
%$f(x)$ agrees with
%\[ (x-.900000000000006)^{18}(x-.99999999999997)^{10}(x-1.10000000000001)^{16}
%\]
%for at least 15 digits.
The condition number is $60.4$. Therefore,
with perturbation at the 16-th digit of the coefficients,
14 correct digits constitute the best possible accuracy that can be expected 
from any method. %Our method achieves this optimal precision.

An important feature of Algorithm I is that it does {\em not} require the 
correct multiplicity structure. 
Computation with different structures is permissible 
with Algorithm I and often needed when the structure is unclear.
If the computation is on a ``wrong'' pejorative manifold, then either the  
condition number or the backward error becomes large. 
Table \ref{table_ex1} is a partial list of pejorative roots under different 
multiplicity structures.
%
%{\bf %%%%
The apparent deviations on the pejorative roots in comparison with
$(0.9,1.0,1.1)$ are the effect of the difference in structure, 
not the failure of the the error bounds in \S \ref{sec:cond} where
structure preservation is assumed.
%} %%%%

\begin{table}[htbp]
\scriptsize
\begin{verbatim}
      multiplicity       pejorative           backward error           condition
       structure           roots               (relative)                number

      [1,1,...,1]    (see Figure 6)          .0000000000000006     1390704851032436
      [18,10,16]   (.9000, 1.0000, 1.1000)   .000000000000002                    60.4
      [17,11,16]   (.8980,  .9934, 1.1006)   .0000004                            53.8
      [14,16,14]   (.8890,  .9892, 1.1090)   .000003                             29.0
      [10,24,10]   (.8711,  .9906, 1.1315)   .000008                             26.7
      [2, 40, 2]   (.7390,  .9917, 1.3277)   .00009                              23.6
      [1, 43]      (.5447, 1.0054)           .004                                 1.3
         [44]             ( .9925)           .04                                   .0058
\end{verbatim}
\normalsize \vspace{-5mm}
\caption{\footnotesize 
Partial list of multiple roots on different pejorative manifolds}
\label{table_ex1} \vspace{-2mm}
\end{table}

%{\bf %%%%
For any polynomial of degree $n$, 
the nearest pejorative manifold(s) {\em always} include
$\Pi_\ell \equiv \bdC^n$ with structure $\ell = [1,1,\cdots,1]$
because every other manifold is its subset. 
For that reason, an {\em unconstraint minimization} of the backward error
(i.e. the distance to a pejorative manifold) will naturally lead to 
the simple, clustered, and incorrect roots as shown
in Figure \ref{trizero}.
Notice that pejorative roots with different multiplicity structures 
corresponds to different backward errors and condition numbers 
$\kappa_{\ell, w}$. 
Generally, reduction in sensitivity may exclude root sets of higher
backward accuracy.
In this example, minimizing the backward error among all pejorative roots
with a sensitivity constraint, say $\kappa_{\ell,w}(\bdz) < 100$, 
leads to the accurate
roots with correct multiplicity structure.
In short, conventional methods seek {\em unconstraint
minimization} of the backward error among all pejorative roots, while 
computing multiple roots of inexact polynomials 
may be accomplished as a {\em constraint optimization} problem that minimizes
the distance to a manifold subject to the condition that the roots are 
insensitive to perturbation with respect to the structure. 
%}%%%%

\subsubsection{Roots with huge multiplicities}

The accuracy as well as stability of Algorithm I seem independent of the
multiplicities of the roots. 
For instance, let's consider the 
polynomial of degree $1000$
\[ g(x) = [x-(0.3+0.6i)]^{100}\,[x-(0.1+0.7i)]^{200}\,
[x-(0.7+0.5i)]^{300}\,[x-(0.3+0.4i)]^{400}. \]
\begin{wrapfigure}[13]{r}{3.0in}
\centerline{ 
\epsfig{figure=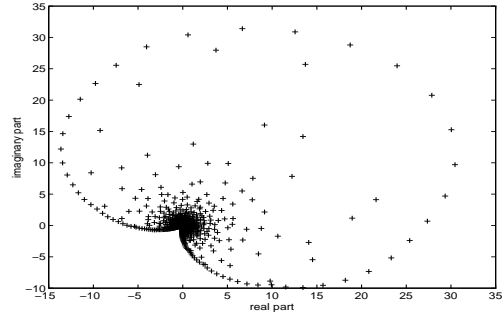,
width=2.6in,height=1.65in}}
\caption[Matlab result for a polynomial of large
degree and high multiplicities]{\footnotesize 
Result for the degree 1000 polynomial by Matlab function {\tt roots}}
\label{fig:ex04}
\end{wrapfigure}
The multiplicities of the roots are $100$, $200$, $300$ and $400$.
These multiplicities are ``huge'' compared to other numerical examples, 
usually with multiplicities less than ten, used in the 
root-finding literature.
In addition to such high multiplicities, we perturb the sixth digits of all 
coefficients of ~$g$~ by multiplying $(1\pm 10^{-6})$ on each one of them.
Using any conventional approach, this perturbation will result in
a total loss of forward accuracy, even if multiprecision is used.
The code {\sc PejRoot} of Algorithm I takes a few seconds under Matlab to 
calculate all roots up to 7 digits accuracy, as iterates
shown in Table \ref{d1000}.
Taking the condition number 0.58 into account, this accuracy
is optimal. On the same machine, Matlab function 
{\tt roots} takes about 15 minutes to produce 1000 incorrect roots, 
as shown in Figure \ref{fig:ex04}.

%\scriptsize
%\begin{tabular}{ll|ll|ll|ll}
%&{\normalsize $z_1$}  && {\normalsize $z_2$} && {\normalsize $z_3$} &&
%{\normalsize $z_4$}  \\ \hline
%.289  & +.601i & .100   &  +.702i & .702 &    +.498i & .301  & +.399i \\
%.309  & +.602i &.097 & +.698i & .698 & +.499i & .299 & +.400i \\
%.293  & +.596i & .101 & +.7003i & .7002 & +.5005i & .3007 & +.4003i\\
%.3003 & +.5994i & .09994 & +.70008i & .69996 &+.50003i &.29996   &+.40007i \\
%.300005 &+.600006  &.099998  &+.6999992i &.69999992 &+.4999993i &
%.2999992 &+.3999992i \\
%.3000002&+.60000005i &.09999995&+.69999998i &.69999997&+.49999998i  &
%.29999997&+.400000002i
%\end{tabular}
%\normalsize

\begin{table}[ht]
\scriptsize
\begin{verbatim}
  .289     +.601i      | .100      +.702i      | .702      +.498i      | .301      +.399i 
  .309     +.602i      | .097      +.698i      | .698      +.499i      | .299      +.400i 
  .293     +.596i      | .101      +.7003i     | .7002     +.5005i     | .3007     +.4003i
  .3003    +.5994i     | .09994    +.70008 i   | .69996    +.50003i    | .29996    +.40007i 
  .300005  +.600006    | .099998   +.6999992i  | .69999992 +.4999993i  | .2999992  +.3999992i 
  .3000002 +.60000005i | .09999995 +.69999998i | .69999997 +.49999998i | .29999997 +.400000002i
\end{verbatim}
\vspace{-4mm}
\caption{\footnotesize Iterates on the degree 1000 polynomial} \label{d1000}
\end{table}

\vspace{-7mm}
\section{Algorithm II: the multiplicity structure and initial root
estimates}
\label{part:gcd}

While Algorithm I can be used on any particular pejorative manifold,
of course, the ``correct'' multiplicity structure is preferred
if it is attainable.
We present Algorithm II that calculates the multiplicity structure of a
given polynomial as well as the initial root approximation for Algorithm I.
 
\subsection{Remarks on the univariate GCD computation}

For a given polynomial $p$ with $u = GCD\left(p,p'\right)$,
%\begin{equation} \label{polyp}
% p(x) = p_0(x-z_1)^{\ell_1}\, (x-z_1)^{\ell_1} \,\cdots \,(x-z_m)^{\ell_m}
%\end{equation}
%with $z_j$'s distinct,
Lagrange pointed out in 1769 that ~$v = p/u$
%~$\displaystyle v = \frac{p}{u}$~
has the same distinct roots as $p$, and all roots of $v$ are simple.
%Namely,
%\[ v(x) = v_0(x-z_1)(x-z_2)\cdots (x-z_m), \;\;\;
%u(x) = u_0(x-z_1)^{\ell_1-1}\cdots (x-z_m)^{\ell_m-1}. \]
If $v$ is obtainable, its simple roots can be calculated
using standard root-finders. 
Based on this observation, the following process that is described
by Gauss in 1863 is 
a natural approach to attain total factorization of the polynomial $p$.
\begin{equation} \label{gfact}
\left\{
\begin{array}{l}
u_0 = p \\
\mbox{for $j=1,2,\cdots $, \ while $deg(u_{j-1}) > 0$ do} \\ %\\
\;\;\;\;\;\;\;\; \mbox{calculate } u_j = GCD\left(u_{j-1}, u_{j-1}'\right), 
\;\;\;\;\; {\displaystyle v_j = \frac{u_{j-1}}{u_j}} \\  %\\
\;\;\;\;\; \mbox{\ \ \ calculate the (simple) roots of } v_j(x) \\
\mbox{end do}
\end{array} \right.
\end{equation}
In this process, 
a $k$-fold root of $p$ will be calculated $k$ times 
as simple roots of $v_j$'s.
Other squarefree factorization processes, such as Yun's algorithm 
\cite{yun}, have also been proposed in the context of Computer Algebra.

The difficulty in carrying out the process (\ref{gfact}) is the GCD 
computation. 
The classical Euclidean GCD Algorithm requires recursive polynomial 
division which may not be numerically stable (see \S \ref{sec:lsqdiv}). 
Therefore, implementations of (\ref{gfact}) based on the Euclidean GCD-finder
\cite{brugnano-trigiante,uhlig-99}
may fail to reach desirable reliability or accuracy
(see numerical comparison in \S 
\ref{sec:nrI}).

Numerical GCD computation has been studied extensively
\cite{chin-corless,corless-gianni,emiris-galligo-lombardi,hribernig-stetter,
karmarkar-lakshman, pan96, rupprecht}.
However, a reliable blackbox-type software is still 
not available.
In \cite{corless-gianni}, Corless, Gianni, Trager and Watt proposed a 
novel approach using the singular value 
decomposition in finding the degree of the GCD, and suggested the
possibility of solving a GCD system similar to (\ref{gcdsys}) below as a 
least 
squares problem, along with several other possibilities including using
the Euclidean algorithm.
%
%In short, Corless-Gianni-Trager-Watt approach consists of two stages. 
%
%First, it uses SVD to determine the rank of the Sylvester resultant 
%matrix (i.e. the $S_{n-1}(p)$ in Definition \ref{def:syl} in our case), 
%thereby obtaining the degree of the GCD. 
%
%Secondly, the coefficients of the GCD are to be determined.
%
%They suggested several possible methods to calculate
%the coefficients of the GCD in the second stage, 
%including the ordinary Euclidean Algorithm 
%using long division, 
%as well as ``standard optimization algorithms'' solving a quadratic 
%least squares problem of a GCD system (e.g. system (\ref{gcdsys}) in our 
%case). 

There are several unresolved issues in the approach of Corless et al,
especially in the stage of calculating the 
GCD after determining its degree. 
Among the possible avenues suggested, they seem to prefer using iterative 
methods to solve the least squares problem similar to (\ref{gcdsys}) below. 
However, their least squares system is underdetermined by one equation.
Moreover, with no clearly decided initial iterate being given, one can only 
leave this
crucial ingredient to guessing or some sort of expensive global 
search \cite{chin-corless}. 
From \cite{corless-gianni} and its follow-up work 
such as \cite{chin-corless,karmarkar-lakshman}
it is also not clear which standard optimization algorithm should be
selected.
We shall demonstrate that the Gauss-Newton iteration, absent from the above
works, is apparently the simplest, most efficient and most suitable method 
in solving the GCD system (\ref{gcdsys}), and it is at least 
locally convergent.

The key to carrying out the procedure (\ref{gfact})
%, thereby calculating 
%roots {\em and} their multiplicities of $p(x)$,
is the capability to factor an arbitrary polynomial $f$ and its 
derivative $f'$ 
with a GCD triplet $(u,v,w)$:
\begin{equation} \label{gcdtriple}
\left\{ \begin{array}{lcl}
u(x)\,v(x) & = &f(x) \\ u(x)\,w(x) &= &f'(x) \end{array} \right.,
\;\;\; \mbox{$u$ is monic, $v$ and $w$ are co-prime}.
\end{equation}
In light of the Corless-Gianni-Trager-Watt approach, which calculates 
all singular values of a single Sylvester matrix $S_{n-1}(f)$, we employ
a successive updating process that calculates only the smallest singular 
values of the Sylvester matrices $S_j(f)$, $j=1,2,\cdots$~ and 
stop at the first rank deficient matrix $S_k(f)$.
With this $S_k(f)$, not only the degrees of the GCD triplet $u,v,w$ 
are available, we also obtain coefficients of $v$ and $w$ automatically
from the resulting right singular vector. 
In combination with a least squares division in \S \ref{ini_it} 
of the unstable long division,
we can generate an approximation to the GCD triplet, 
and obtain an initial iterate that is not clearly indicated in 
the approach of Corless et al. 
Consequently, a blackbox-type software computing $GCD(f,f')$
is developed for the process (\ref{gfact}).

%The GCD computation is an essential component of our Algorithm II. 
%
The discussion on GCD computation is limited to $GCD(f,f')$ 
because our objective is mainly root-finding in this paper.
With minor modifications, our GCD-finder can easily be adapted to 
the general GCD problem of arbitrary polynomial pairs.

\subsection{Calculating the greatest common divisor} \label{sec:gcd}

Algorithm II is based on the following  GCD-finder for an arbitrary polynomial
$f$:

\vspace{-8mm}
\begin{tabbing}
\hspace{3mm} \=
\hspace{6mm} \=
\hspace{6mm} \=
\hspace{6mm} \=
\hspace{6mm} \=
\hspace{6mm} \=
\hspace{6mm} \= \\
\> {\tt \bf STEP 1.} Find the degree $m$ of $GCD(f,f')$. \\
\> {\tt \bf STEP 2.} 
Set up the system (\ref{gcdtriple}) in accordance with the degree $m$. \\
%for the GCD triplet $(u, v, w)$. \\
\>{\tt \bf STEP 3.} 
Find an initial approximation to $u$, $v$ and $w$ for the GCD system
(\ref{gcdtriple}). \\
\> {\tt \bf STEP 4.} 
Use the Gauss-Newton iteration to refine the GCD triplet $(u, v, w)$.
\end{tabbing}

\vspace{-4mm}
We shall describe each step in detail.

\subsubsection{Finding the degrees of the GCD triplet} \label{sec:deg}

Let $f$ be a polynomial of degree $n$. 
By Lemma \ref{lem:gcdsv}, the degree of $u = GCD(f,f')$ is $m=n-k$ 
if and only if the the $k$-th Sylvester discriminant matrix is 
the first one being rank-deficient.
Therefore, $m = deg(u)$ can be identified by calculating the sequence of the 
smallest singular values $\vs_j$ of $S_j(f)$, $j=1,2,\cdots$, until reaching
$\vs_k$ that is {\em approximately} zero. 
Since only one singular pair (i.e. the singular value and the right singular 
vector) is needed, the inverse iteration described 
in Lemma \ref{lem:invit} is suitable for this purpose.
Moreover, we can reduce the computing cost even further by recycling and 
updating the QR decomposition of $S_j(f)$'s along the way. 
More specifically, let 
\[ f(x) = a_0 x^n + a_1 x^{n-1} + \cdots + a_n, \;\;\; 
f'(x) = b_0 x^{n-1} + b_{1} x^{n-2} + \cdots + b_{n-1}. \]
We rotate the columns of $S_j(f)$ to form $\hat{S}_j(f)$ in such a way
\footnotesize
\begin{eqnarray*}
\lefteqn{
\begin{array}{c}
\;\;\; \overbrace{\mbox{\ \ \ \ \ \ \ \ \ \ \ \ \ \ 
\ \ \ \ \ \ \ }}^{j+1} \;\;\;\;\;\;
\overbrace{\mbox{\ \ \ \ \ \ \ \ \ \ \ \ 
\ \ \ \ \ \ \ }}^{j} \\
\left( \begin{array}{llllll}
b_0  & & & a_0 && \\
b_1 & \ddots  & & a_1 & \ddots & \\
\vdots & \ddots & b_0 & \vdots & \ddots & a_0 \\
b_{n-1} & & b_1 & a_n  & & a_1 \\
 & \ddots & \vdots &  & \ddots & \vdots \\
&        & b_{n-1} &&        & a_n   
\end{array} \right)
\end{array} } \\
&& \longrightarrow  \;\;\;\;
\begin{array}{c}
\;\;\; \overbrace{\mbox{\ \ \ \ \ \ \ \ \ \ \ \ \ \ \ \ \ \ \ \ \ \ 
\ \ \ \ \ \ \ \ \ \ \ \ \ \ 
\ \ \ \ \ \ \ \ \ \ \ 
\ \ \ \ \ \ \ \ \ \ \ \ \ \ \ \ \ \ \ \ \ \ \ \ \ \ }}^{2j+1}  \\
\left( \begin{array}{llllllll} 
b_0 & a_0 &&&&&& \\
b_1 & a_1 & b_0 & a_0 &&&& \\
\vdots & \vdots & b_1 & a_1 & \ddots &&& \\
b_{n-1} & a_{n-1} & \vdots & \vdots & & b_0 & a_0 & \\
& a_n & b_{n-1} & a_{n-1} & & b_1 & a_1 & b_0 \\
& & & a_n  & \ddots & \vdots & \vdots & b_1 \\
&&&&& b_{n-1} & a_{n-1} & \vdots \\
&&&&&& a_n & b_{n-1} 
\end{array} \right)
\end{array}
\end{eqnarray*}  \normalsize
that the odd and even columns of $\hat{S}_j(f)$ consist of the coefficients of 
$f'$ and $f$ respectively. 
Consequently, the matrix $\hat{S}_{j+1}(f)$ is simply formed by adding 
a zero row at the bottom and two columns in the right on 
$\hat{S}_{j}(f)$. 
Updating the QR decomposition of each $\hat{S}_j(f)$ 
requires only $O(n)$ additional flops. 
The inverse iteration (\ref{invit}) 
%consisting of a forward and a backward substitutions 
requires $O(j^2)$ flops at each $S_j(f)$.

Let $\theta$ be a given {\em zero singular value threshold},
which shall be discussed more in \S \ref{cpara}.
With successive QR updating and the inverse iteration, 
the process of finding the degrees of the GCD triplet $(u,v,w)$ 
can be summarized as follows.  

\vspace{-3mm}
\begin{tabbing}
\hspace{0mm} \=
\hspace{6mm} \=
\hspace{6mm} \=
\hspace{6mm} \=
\hspace{6mm} \=
\hspace{6mm} \=
\hspace{6mm} \= \\
\>\> Calculate the QR decomposition of the $(n+1)\times 3$ matrix 
$\hat{S}_1(f) = Q_1R_1$
\\
\>\> {\tt For} ~$j = 1, 2, \cdots$~ {\tt do} \\
\>\>\> use the inverse iteration (\ref{invit}) to find 
the smallest singular value $\vs_j$ \\
\>\>\>\> of $\hat{S}_j(f)$ and the corresponding 
right singular vector $\bdy_j$\\
\>\>\> {\tt if} ~$\vs_j \leq \theta \| \,\bdf \,\|_2$~ {\tt then} 
%\>\>\>\> 
~~~$k = j$, \  $m=n-k$, \ extract $\bdv$ and $\bdw$ from $\bdy_j$, 
\ {\tt exit} \\
\>\>\> {\tt else} %\\
%\>\>\>\> 
\ \ \ \  update ~$\hat{S}_j(f)$~ to ~$\hat{S}_{j+1}(f) = Q_{j+1}R_{j+1}$ \\
\>\>\> {\tt end if} \\
\>\> {\tt end do} 
\end{tabbing}

\subsubsection{The quadratic GCD system} \label{sec:gcdsys}

Let $m = n-k$ be the degree of $GCD(f,f')$ calculated in {\tt STEP 1}.
We now formulate the GCD system (\ref{gcdtriple}) of {\tt STEP 2}
in vector form with unknown vectors
$\bdu$, $\bdv$ and $\bdw$:
\begin{equation} \label{gcdsys}
\left[ \begin{array}{l} u_0 \\ conv(\bdu,\bdv) \\ conv(\bdu,\bdw)
\end{array} \right]
 = \left[ \begin{array}{l} 1 \\  \bdf \\ \bdf'
\end{array} \right], \;\;\;\; \mbox{ \ for \ }
\left(\begin{array}{c} \bdu \\ \bdv \\ \bdw \end{array} \right) 
\in \bdC^{m+1} \times \bdC^{k+1} \times \bdC^{k}.
\end{equation}
Here, the convolution $conv(\cdot,\cdot)$ is defined in 
Lemma \ref{lem:conv}.
The following lemma ensures this quadratic system is
nonsingular.

\begin{lemma} \label{lem:jfr}
~The Jacobian of the quadratic system
%$F(\bdu,\bdv,\bdw) = 
%\left[ \begin{array}{l} u_0 \\ conv(\bdu,\bdv) \\ conv(\bdu,\bdw)
%\end{array} \right] $ 
(\ref{gcdsys}) is
\begin{eqnarray} \label{syljac}
 J(\bdu,\bdv,\bdw) & = & \left[ 
\begin{array}{cll} 
\bde_1^\top &  & \\
C_{m}(v) & C_{k}(u) &  \\
C_{m}(w) &  & C_{k-1}(u)
%\bde_1^\top & 0_{1 \times (k+1)} & 0_{1\times k} \\
%C_{m}(v) & C_{k}(u) & 0_{(n+1)\times k} \\
%C_{m}(w) & 0_{n\times (k+1)}  & C_{k-1}(u)
\end{array}
\right], \\
\mbox{\ \ where } && \bde_1 = (1,0,\cdots,0)^\top \in \bdC^{m+1}.
\nonumber
\end{eqnarray}
If ~$u = GCD(f,f')$ with $(\bdu,\bdv,\bdw)$ satisfying (\ref{gcdsys}),
then ~$J(\bdu,\bdv,\bdw)$ is of full (column) rank.
\end{lemma}

{\em Proof.} It is straightforward to verify (\ref{syljac}) by using Lemma
\ref{lem:conv}. 
 To prove $J(\bdu,\bdv,\bdw)$ is of full rank, we assume the
existence of polynomials $q(x) = \sum_{j=0}^m q_j x^{m-j}$, 
~$r(x) = \sum_{j=0}^k r_j x^{k-j}$ and 
$s(x) = \sum_{j=0}^{k-1} s_j x^{k-j-1}$~
such that 
\begin{equation} \label{Jqrs}
 J(\bdu,\bdv,\bdw) \left( \begin{array}{c} \bdq \\ \bdr \\ \bds 
\end{array} \right) = 0, \mbox{ \ \ or \ \ } 
\left\{ \begin{array}{lcc} q_0 & = & 0 \\
vq + ur & = & 0 \\ wq + us & = & 0 \end{array} \right. .
\end{equation}
Here, as before, $\bdq$, $\bdr$ and $\bds$ are coefficient vectors of $q$, $r$ 
and $s$ respectively. 
From (\ref{Jqrs}), we have $vq = -ur$ and $wq=-us$.
So, $wvq-vwq = -uwr + uvs = 0$, namely $-wr + vs = 0$
or $wr=vs$.
Since $v$ and $w$ are co-prime, 
there is a polynomial $t$ such that $r=tv$ and $s=tw$.
Consequently, $vq = -ur = -utv$ leads to $q = -tu$. Because 
$deg(q) = deg(tu) \leq m$,
$deg(u)=m \geq 0$ and $u_0=1$, the degree of $t$ must be zero. 
So the polynomial $t$ is a constant.
Using the first equation in (\ref{Jqrs}) and $u_0=1$,  we have
$q_0 = -tu_0 = -t = 0$. It follows that 
$q=-tu=0$, $r=tv=0$ and $s=tw =0$. 
Consequently, $J(\bdu,\bdv,\bdw)$ is of full rank.
\qed

%{\bf %%%%
The equation $u_0=1$, 
absent in Corless-Gianni-Trager-Watt GCD method \cite{corless-gianni}, 
needs to be included either explicitly 
with an equation or implicitly by eliminating
$u_0$ from the variables to ensure the regularity of
the system (\ref{gcdsys}). 
In Remark 1 of \cite[\S 4.3]{chin-corless} with a heuristic explanation, 
Chin, Corless and Corless realized that the restriction $u_0=1$ 
may make their Divisor-Quotient Iteration converge, but abandoned it 
since their ``test showed that the overall performance was worse when 
this constraint was in place'' \cite[\S 5.1.4]{chin-corless}. 
Because we use a different refinement approach in our GCD-finder, 
preserving this constraint, 
and thereby the regularity of the GCD system
(\ref{gcdsys}), may be the very reason for our method to obtain more 
robust test results.
Without this regularity, 
the local convergence of the Gauss-Newton iteration we use would not be
guaranteed. 
%} %%%%

\begin{theorem} \label{gnconv2}
~Let ~$\tilde{u} = GCD(f,f')$~ 
with ~$\tilde{v}$~ and ~$\tilde{w}$~ satisfying (\ref{gcdsys}),
and let ~$W$~ be a weight matrix.
Then there exists ~$\eps>0$~ such that for all ~$\bdu_0$, $\bdv_0$, $\bdw_0$
satisfying $\bnorm{\bdu_0-\tilde{\bdu}} < \eps$, 
$ \bnorm{\bdv_0-\tilde{\bdv}} < \eps$~ and
$\bnorm{\bdw_0-\tilde{\bdw}} < \eps$, 
the Gauss-Newton iteration 
\begin{eqnarray} \label{G-N_it1}
\;\;\;\;
\left[ \begin{array}{c} \bdu_{j+1} \\ \bdv_{j+1} \\ \bdw_{j+1} 
\end{array} \right] & = &
\left[ \begin{array}{c} \bdu_{j} \\ \bdv_{j} \\ \bdw_{j} 
\end{array} \right]
- J(\bdu_j,\bdv_j,\bdw_j)_W^+ \left[\,
 \begin{array}{lcl} 
\bde_1^\top \bdu_j & - & 1 \\
conv(\bdu_j,\bdv_j) & - & \bdf \\ 
conv(\bdu_j,\bdw_j) & - & \bdf'
\end{array} 
 \,\right] \\ &&  \nonumber \\
 & & \;\;\; j = 0, 1, \cdots \nonumber
\end{eqnarray}
converges to $[\tilde{\bdu},\tilde{\bdv},\tilde{\bdw}]^\top$ quadratically.
Here $J(\cdot)^+_W = [J(\cdot)^H W^2 J(\cdot)]^{-1} J(\cdot)^H W^2$ is the
weighted pseudo-inverse of the Jacobian $J(\cdot)$ as
defined in (\ref{syljac}).
\end{theorem}

{\em Proof.} 
A straightforward verification by using Lemma \ref{conv_lem} and Lemma 
\ref{lem:jfr}.
\qed

\subsubsection{Setting up the initial iterate} \label{ini_it}
\label{sec:lsqdiv}

We now need initial iterates ~$\bdu_0,\bdv_0,\bdw_0$ for the 
Gauss-Newton iteration (\ref{G-N_it1}).
In {\tt STEP 1}, when the singular value $\vs_k$ is calculated,
the associated singular vector $\bdy_k$ 
consists of $\bdv_0$ and $\bdw_0$ that are approximations to 
$\bdv$ and $\bdw$ in (\ref{gcdsys}) respectively (see Lemma \ref{lem:uvw}). 
%
%For the purpose of refinement, we shall use these vectors as our initial 
%approximation $\bdv_0$ and $\bdw_0$. 
%
Because of the column rotation in \S \ref{sec:deg}, 
the odd and even entries of $\bdy_k$ form $\bdv_0$ and $\bdw_0$ respectively.
For the initial approximation $\bdu_0$,
notice that in theory the long division yields, 
\begin{equation} \label{ldiv0}  f(x) = v_0(x)q(x) + r(x) \end{equation}
with $u_0(x) = q(x)$ and $r(x) = 0$.
The process itself may not be numerically stable. 

%{\bf %%%%
In the context of Corless-Gianni-Trager-Watt method, 
Corless et al. \cite[Lemma 3]{corless-gianni} propose the use of
least squares method to minimize $\|\Delta \bdf\|_2 = 
\|\,C_{n-m}(d) \bdh - \bdf \, \|_2$ whenever a candidate $d$ of 
degree $m$ approximating $GCD(f,g)$ is available. 
In our approach, there is no candidate for $u=GCD(f,f')$ at the end of 
{\tt STEP 1}. Instead, we have $v_0$ and $w_0$ available 
that approximate $v$ and $w$ respectively such that $u = f/v = f'/w$. 
Whereas we can adapt the least squares strategy, which Corless et al.
use to calculate $\|\Delta \bdf\|_2$ from given $u_0$,
to calculate the approximation $u_0$ of $GCD(f,f')$ 
from given $v_0$ and $w_0$ in our method, and justify it via 
a condition theory of linear systems.
%} %%%%

By Lemma \ref{lem:uvw}, the long division (\ref{ldiv0}) with $r(x)=0$
is equivalent to solving the linear system
\begin{equation} \label{lsqdiv} C_{m}(v_0) \, \bdu_0 = \bdf \end{equation}
for a least squares solution $\bdu_0$
that minimizes $\bnorm{conv(\bdu_0,\bdv_0)-\bdf}$.
This ``least squares division'' is more accurate than the long division
(\ref{ldiv0}).
In fact, the long division (\ref{ldiv0}) is equivalent to solving the
$(n+1)\times (n+1)$ lower triangular linear system
\begin{equation} \label{ldivsys}
L_{m}(v_0)
\left( \begin{array}{c} \bdq \\ \bdr \end{array} \right) = \bdf,
\mbox{ \ \ with \ \ } L_{m}(v_0) =
\left( \begin{array}{cc}
C_{m}(v_0) & \left| \begin{array}{c} 0_{(m+1)\times (n-m)} \\
I_{(n-m)\times (n-m)}
\end{array} \right.
\end{array} \right).
\end{equation}
The following theorem indicates that
solving (\ref{lsqdiv}) for $\bdu_0$ may be more preferable than using
the long division (\ref{ldiv0}).

\begin{theorem}
Let $\kappa(A)$ denote the condition number of an arbitrary matrix $A$ with
respect to the matrix 2-norm. Then
~~\( \kappa\left(C_{m}(v)\right) \leq \kappa(L_{m}(v)) \)~
for any polynomial $v$ and $m>0$.
\end{theorem}

{\em Proof.}
~For any matrix $A$,
~$\kappa(A) = \frac{\sg_{max}(A)}{\sg_{min}(A)}$, where
\[ \displaystyle \sg_{max}(A) = \max_{\|\bdx\|_2 = 1} \bnorm{A\bdx}
\mbox{ \ \ and \ \ }
\sg_{min}(A) = \min_{\|\bdx\|_2 = 1} \bnorm{A\bdx} \]
are the largest and smallest singular values of $A$ respectively.
Therefore
\begin{eqnarray*}
\lefteqn{
\sg_{max}(C_{m}(v))  =  \max_{\|\bdu\|_2 = 1} \bnorm{C_{m}(v)\,\bdu}
= \max_{\|\bdq\|_2 = 1,\bdr=0} \bnorm{C_{m}(v)\,\bdq + \bdr}} \\
& = & \max_{\|\bdq\|_2 = 1,\bdr=0} \left\| L_{m}(v)
\left( \begin{array}{c} \bdq \\ \bdr \end{array} \right) \right\|_2
\;\; \leq \;\; \max_{\|y\|_2 = 1} \bnorm{L_{m}(v) \bdy} =
\sg_{max} (L_{m}(v)).
\end{eqnarray*}
Similarly, ~$\sg_{min}(C_{m}(v)) \geq \sg_{min}(L_{m}(v))$, and
consequently, $\kappa(C_{m}(v)) \leq \kappa(L_{m}(v))$.
\qed

The magnitude gap between the condition numbers $\kappa(C_{m}(v))$ and
$\kappa(L_{m}(v))$ can be tremendous for seemingly harmless $v$
and moderate $m$.
Actually, $L_{m}(v)$ can be pathetically ill-conditioned, 
making the long division (\ref{ldiv0}) virtually a singular process, 
while $C_{m}(v)$ is still well conditioned. 
For example, consider a simple polynomial $v(x) = x + 25$. 
When $m$ increases, $\kappa(L_{m}(v))$ grows exponentially but 
$\kappa(C_{m}(v))$ stays as nearly a constant, see Table \ref{divcond}.
\begin{table}[ht] 
\begin{center}
\begin{tabular}{|c||r|l|l|l|} \hline
 & $m=1$ & $m = 5$ & $m = 10$ & $m=20$ \\ \hline
$\kappa(C_{m}(v))$ & 1 & 1.0668 & 1.0791 & 1.0823  \\
$\kappa(L_{m}(v))$ & 627 & $1.01 \times 10^7$ &  $9.92 \times 10^{13}$ & 
$9.46 \times 10^{27}$ \\ \hline
\end{tabular} 
\end{center}  \vspace{-4mm}
\caption{\footnotesize
The comparison between the conditions of (\ref{ldiv0}) and
(\ref{lsqdiv}) for $v(x) = x+25$.}\label{divcond}
\vspace{-3mm}
\end{table}
In fact, we have not encountered a truly ill-conditioned least squares 
division (\ref{lsqdiv}) in our extensive numerical experiments.
On the other hand, the example shown in Table \ref{ldiv} is quite common.
In which ~$\bdf = conv(\bdu,\bdv)$~ is rounded up at the eighth digit
after decimal point.
The difference between the long division (Matlab {\tt deconv}) 
and the least squares division is quite substantial.

\begin{table}[ht]
\scriptsize
\begin{center}
\begin{tabular}{||r|r||c|c|l||} \hline
\multicolumn{2}{||c||}{Data} & 
\multicolumn{3}{c||}{Comparison} \\ \hline
approx. coef. & coefficients \ \  & known coef.'s of \ \ & least squares 
& \ \ \ \ long \\
of $f(x)$ \ \ \ & of $v(x)$ \ \ \ & $f(x)\div v(x)$ & division \ \  
& \ \ division \\ 
\hline
   1.00000000 &  1.00000000 & 1.00000000 & 0.9999999999 & \ \ 1.00000000 \\
  23.35360257 & 23.01829201 & 0.33531056 & 0.3353105599 & \ \ 0.33531056 \\
  29.89831582 & 22.05776405 & 0.12227539 & 0.1222753902 & \ \ 0.122275385 \\
  10.75803809 &             & 0.54726624 & 0.5472662398 & \ \ 0.5472663 \\
  15.57240922 &             & 0.27815340 & 0.2781534002 & \ \ 0.278151 \\
  18.76038493 &             & 0.28629915 & 0.2862991496 & \ \ 0.28634 \\
  13.73079603 &             & 1.00523653 & 1.0052365305 & \ \ 1.004 \\
  30.45600101 &             & 1.00205392 & 1.0020539195 & \ \ 1.02 \\
  46.21275197 &             & 0.97391204 & 0.9739120403 & \ \ 0.5 \\
  44.89871211 &             & 0.37785145 & 0.3778514500 & 11. \\
  30.17981700 &&&& \\
   8.33455813 &&&&  \\ \hline
\end{tabular}
\end{center} \vspace{-4mm}
\caption{\footnotesize
A numerical comparison between long division and least squares
division}
\label{ldiv} \vspace{-4mm}
\end{table}

Extracting $\bdv_0$ and $\bdw_0$ from the singular vector and solving
(\ref{lsqdiv}) for $\bdu_0$, we shall use them as the initial iterates for 
the Gauss-Newton iteration (\ref{G-N_it1}) that refines the GCD triplet. 
Moreover, the linear system (\ref{lsqdiv}) is banded, 
with bandwidth being one plus the number of distinct roots. 
Therefore, the cost of solving (\ref{lsqdiv}) is insignificant 
%(no more than $O(n^3))$ 
in the overall complexity. 

\subsubsection{Refining the GCD with the Gauss-Newton iteration}
\label{refgcd} 

The Gauss-Newton iteration is expected to reduce the 
residual 
\begin{equation} \label{gnres}
\left\|\, \left( \begin{array}{l} conv(\bdu_j,\bdv_j) \\ conv(\bdu_j,\bdw_j)
\end{array} \right) - 
\left( \begin{array}{l} \bdf \\ \bdf' \end{array} \right)
\right\|_W 
= \left\|\, W \left( \begin{array}{lcl} conv(\bdu_j,\bdv_j) & - & \bdf \\ 
conv(\bdu_j,\bdw_j) & - & \bdf' \end{array} \right)
\right\|_2
\end{equation}
at each step until it is numerically unreducible. We stop
the iteration when this residual no longer decreases.
The diagonal weight matrix $W$ is used to scale the GCD system 
(\ref{gcdsys}) so that the entries of 
\( W \left[ \begin{array}{l} \bdf \\ \bdf'
\end{array} \right] 
\)
are of similar magnitude.
Each step of the Gauss-Newton iteration 
requires solving an overdetermined linear system
\[
\Big[\; WJ(\bdu_j,\bdv_j,\bdw_j)\;\Big]\;\bdz
%\left( 
%\left[ \begin{array}{c} \bdu_{j+1} \\ \bdv_{j+1} \\ \bdw_{j+1}
%\end{array} \right]
%- \left[ \begin{array}{c} \bdu_{j} \\ \bdv_{j} \\ \bdw_{j}
%\end{array} \right] \right) 
= W \left[\,
\begin{array}{lcl} \bde_1^\top \bdu_j  & - & 1 \\
conv(\bdu_j,\bdv_j) & - & \bdf \\
conv(\bdu_j,\bdw_j) & - & \bdf' \end{array}
\,\right]
\]
\begin{wrapfigure}[11]{r}{2.8in}
\vspace{-3mm}
\centerline{ %\raisebox{1in}{$WJ(\bdu,\bdv,\bdw) \;\;\; = \;\;\; $}
\epsfig{figure=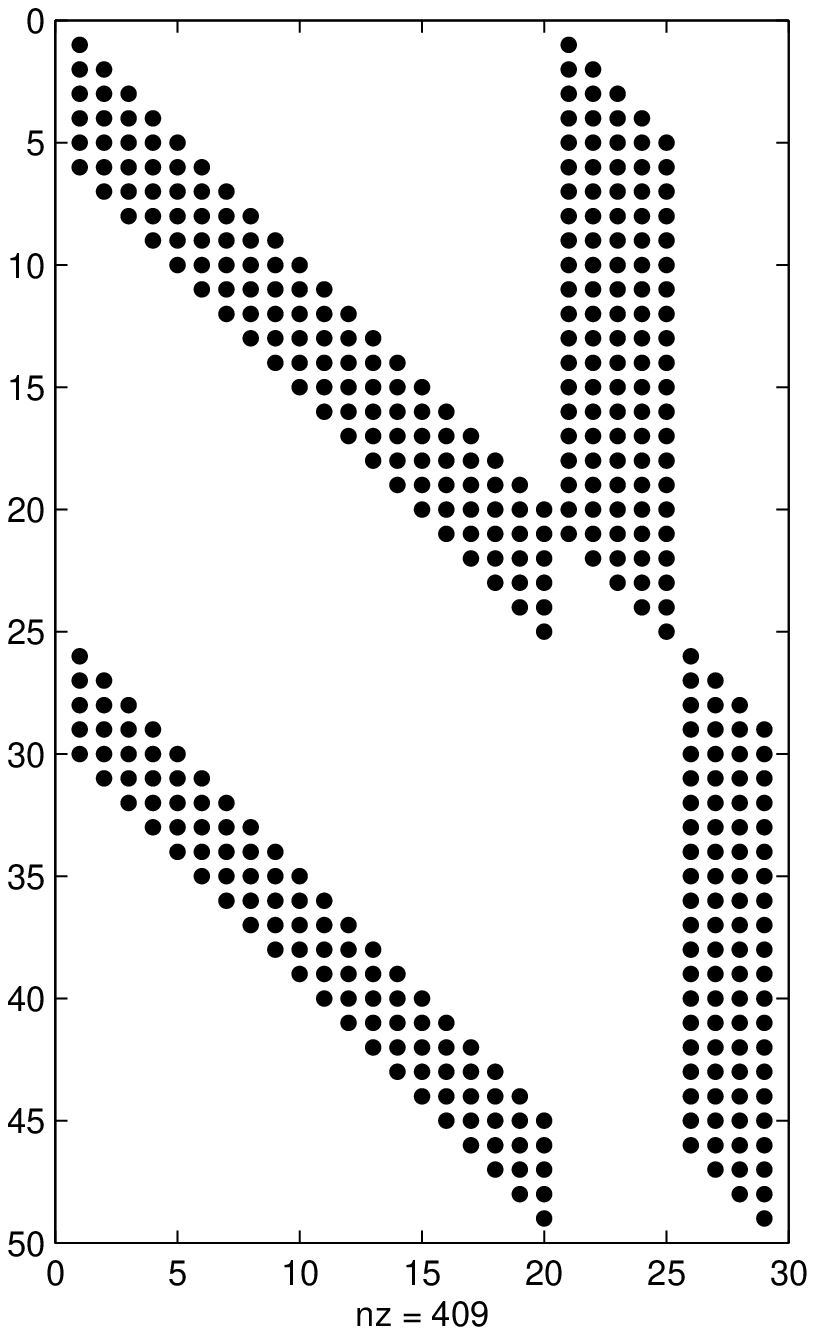,height=1.65in,width=0.8in}
\ \  \raisebox{0.8in}{$\stackrel{QR}{\mbox{\huge $\longrightarrow$}}$} \ \
\epsfig{figure=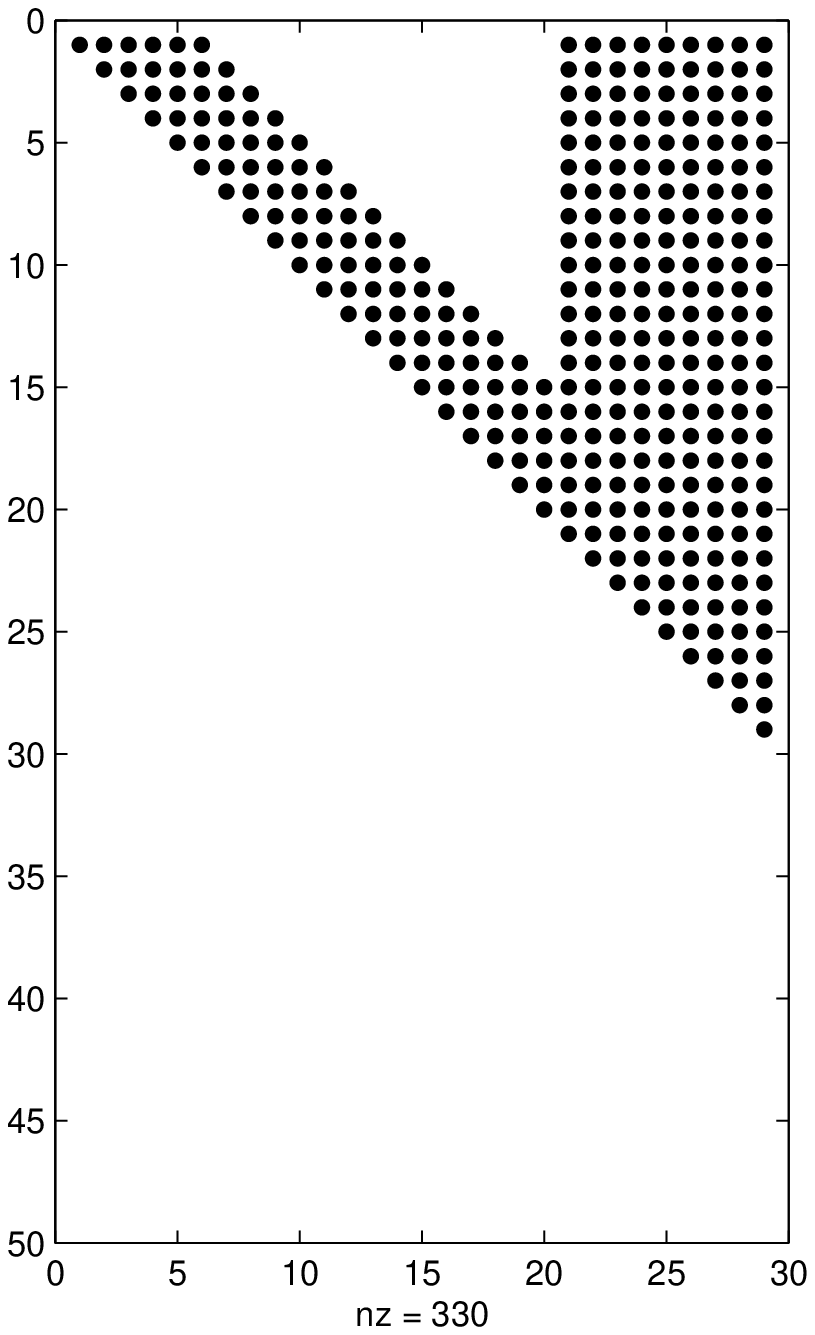,height=1.65in,width=0.8in}
}\vspace{-3mm}
\caption[Sparsity of the Jacobian]{
\footnotesize Sparsity of $J(\bdu,\bdv,\bdw)$
and its triangularization
}
\label{spyjr}
\end{wrapfigure}
for its least squares solution ~$\bdz$, and requires a QR 
decomposition of the Jacobian $WJ(\bdu_j,\bdv_j,\bdw_j)$ and 
a backward substitution for an upper triangular linear system.
This Jacobian is a sparse matrix with a special sparsity structure
that can largely be preserved during the process.
Figure \ref{spyjr} shows the typical sparsity of $WJ(\bdu,\bdv,\bdw)$
along with its triangularization.
When $f$ is a polynomial of degree $n$, a straightforward QR 
decomposition of $WJ(\bdu,\bdv,\bdw)$ costs $O(n^3)$ flops. 
Taking the sparsity of $WJ(\bdu,\bdv,\bdw)$ into account, 
it can be verified that the sparse QR decomposition costs 
$O(mk^2+m^2k+k^3)$, where, as before, $k$ is the number of distinct roots 
and $m=n-k$. 
If $k$ is $o(n)$, then the complexity is reduced to $O(kn^2)$.
%The complexity is significantly reduced to between $O(n^2)$ and $O(n^3)$. 

\subsection{Computing the multiplicity structure and initial root
estimates}
\label{sec:ovI}

%{\bf %%%%
For a given polynomial $p$, the procedure (\ref{gfact}) 
generates a sequence of square-free polynomials
~$v_1, v_2, \cdots, v_s$~ of degrees 
~$d_1 \geq d_2 \geq \cdots \geq d_s$~ respectively such that
~$p = v_1 \; v_2 \; \cdots \; v_s $~ with
\[ \{\mbox{roots of }v_1\} \supseteq \{\mbox{roots of }v_2\} 
\supseteq \cdots \supseteq \{\mbox{roots of }v_s\} \]
Moreover, for each $j = 1, \cdots, s$, all roots of $v_j$ are simple.
Roots of $v_1$ consists of all distinct roots of $p$, while roots of $v_2$ 
consists of all distinct roots of $p/v_1$, etc.
With these properties, the multiplicity structure is determined by
the degrees $d_1,d_2,\cdots,d_s$. For example, consider
\[ p(x) = (x-a)(x-b)^3 (x-c)^4 \]
for any $a$, $b$ and $c$. We have
\begin{center}
\begin{tabular}{crcr}
& $v_j$'s \ \ \ \ \ \  & \ \ \  degrees of $v_j$'s \ \ \  & roots \\ \hline
$v_1(x) = $ & $(x-a)(x-b)(x-c)$ & $d_1\;=\;3$ & $a$ \ \ $b$ \ \ $c$ \\
$v_2(x) = $ & $(x-b)(x-c)$      & $d_2\;=\;2$ & $b$ \ \ $c$ \\
$v_3(x) = $ & $(x-b)(x-c)$      & $d_3\;=\;2$ & $b$ \ \ $c$ \\
$v_4(x) = $ & $(x-c)$           & $d_4\;=\;1$ & $c$ \\ \hline
& multiplicity structure of $p$: &  & 1 \ \ 3 \ \ 4
\end{tabular}
\end{center}
Without locating the roots $a$, $b$ and $c$, the multiplicity structure 
$[\ell_1,\ell_2,\ell_3] = [1,3,4]$ is determined solely from the degrees
$d_1,\cdots,d_4$: 
\begin{eqnarray*}
\ell_1 = 1 &&  \mbox{since \ \ } 
d_1 \;\;\;\;\;\;\;\;\;\;\;\;\;\;\;\;\; \geq 3  
\;\;\;\;\; = (d_1+1)-1\\
\ell_2 = 3 &&  \mbox{since \ \ } d_1, d_2, d_3 \;\;\;\;\;\;\; \geq 2 
\;\;\;\;\; = (d_1+1)-2\\
\ell_3 = 4 &&  \mbox{since \ \ } d_1, d_2, d_3, d_4 \;\; \geq 1 
\;\;\;\;\; = (d_1+1)-3\\
\end{eqnarray*}
Generally, we have the following theorem on identifying the multiplicity 
structure.

\begin{theorem}
For a given polynomial $p$, 
let $v_1,\cdots,v_s$ be the squarefree factors of $p$ generated 
by the procedure (\ref{gfact}) with degrees $d_1 \geq d_2 \geq 
\cdots \geq d_s$ respectively.
Let $k = d_1 = deg(v_1)$. 
Then the multiplicity structure $\ell$ 
of $p$ consists of components 
\begin{equation} \label{strct} 
\ell_j = \max \left\{\,t\,\Big|\, d_t \geq (d_1+1)-j\,\right\}, \;\;\; 
j = 1, 2, \cdots, k.
\end{equation}
\end{theorem}

{\em Proof.} A straightforward verification. \qed

The location of the roots is not needed in deciding the structure.

The initial root approximation is determined based on the fact that
an $l$-fold root of $p(x)$ appears $l$ times as a simple 
root of each polynomial among $v_1, \cdots,v_l$. 
After calculating roots of each $v_j$ with a standard root-finder,
numerically ``identical'' roots of $v_j$'s are grouped in a 
straightforward manner,
according to the multiplicity structure $[\ell_1,\cdots,\ell_k]$
determined by (\ref{strct}), 
to form the initial root approximation $(z_1,\cdots,z_k)$ that is 
needed by Algorithm I.
%} %%%%

\subsection{Control parameters}
\label{cpara}

We use three control parameters for the recursive GCD computation. 
The default values of those parameters given below are selected under the 
assumption that the IEEE standard double precision of 16 decimal digits
is used.
The first control parameter is the {\em zero singular value threshold} 
~$\theta$ for identifying the zero singular value. 
The default choice is $\theta = 10^{-8}$. 
When the smallest singular value 
$\vs_l$ of $\hat{S}_l(u_{m-1})$ is less than $\theta \bnorm{\bdu_{m-1}}$, it
will be {\em tentatively} considered as a zero (pending confirmation
from the residual information produced by the Gauss-Newton iteration). 
Then the Gauss-Newton
iteration is initiated to further reduce the residual as in (\ref{gnres})
to its numerical limit.
We use the second control parameter, the {\em initial residual tolerance}
$\varrho$, to decide if the refined residual is acceptable. 
Our default choice is $\varrho = 10^{-10}$. 
We accept the GCD triplet $(u_m,v_m,w_m)$ when the residual 
\begin{equation} \label{gnres_m}
\rho_m = \left\|\, \left(
\begin{array}{l} conv(\bdu_m,\bdv_m) - \bdu_{m-1} \\
conv(\bdu_m,\bdw_m) - \bdu_{m-1}' \end{array} \right)
\right\|_W \;\; \leq \;\; \varrho \, \bnorm{\bdu_{m-1}}.
\end{equation}
Otherwise, we continue to update $S_l(u_{m-1})$ to $S_{l+1}(u_{m-1})$ and
check $\vs_{l+1}, \cdots$. 

The third parameter is the {\em residual tolerance growth factor} $\phi$. 
Whenever a GCD triplet $(u_m,v_m,w_m)$ and $\rho_m$ are calculated, 
The error in $(u_m,v_m,w_m)$ may cause the residual $\rho_{m+1}$ of 
$(u_{m+1},v_{m+1},w_{m+1})$ to grow. 
Therefore, the tolerance $\varrho$ may need adjustment. 
Our default growth factor is $100$. 
After obtaining $\rho_m$, the residual tolerance $\varrho$ is adjusted 
to be
~\( \max \Big\{\; \varrho, \;\;\phi\, \rho_m \; \Big\}. \)
~Notice that the growth factor is applied to residual $\rho_m$ rather
than the residual tolerance $\varrho$. % to provide some breathing room. 
The residual tolerance $\varrho$ itself may not grow at every step. 

From our computing experience, 
the default control parameters works well for ``normal'' polynomials,
such as those with unclustered roots of moderate multiplicities.
For difficult problems, one may manually adjust the parameters.
The overall Algorithm II shown in Fig. \ref{zgcdroot}~ is implemented
as a Matlab code {\sc GcdRoot} and included in the {\sc MultRoot} package.

\begin{figure}[p] 
\small
\begin{center}
\shadowbox{\parbox{2.0in}{
\baselineskip0mm
\vspace{-0mm}
{\tt
\begin{tabbing}
\hspace{0mm} \=
\hspace{6mm} \=
\hspace{6mm} \=
\hspace{6mm} \=
\hspace{6mm} \=
\hspace{6mm} \=
\hspace{6mm} \= \\
\>{\bf Pseudo-code } {\sc GcdRoot} ~(Algorithm II)\\
\> ~input: The polynomial $p$ of degree $n$, singular threshold $\theta$, \\
\>\>\>\> residual tolerance $\varrho$, residual growth factor $\phi$. \\
\>\>\>\> (If only $p$ is provided, set $\theta=10^{-8}$, $\varrho=10^{-10}$,
            $\phi = 100$ ) \\
\> ~output: the root estimates $(z_1,\cdots,z_k)^\top$ and  \\
\>\>\>multiplicity 
structure $[\ell_1,\cdots,\ell_k]$~ \\
\>\> \\
\>\> Initialize ~$u_0 = p$ \\
\>\> for $m = 1, 2, \cdots, s$, until $deg(u_s) = 0$ ~do \\
\>\>\> for $l = 1, 2, \cdots$ until residual $\rho < \varrho 
\bnorm{\bdu_{m-1}}$ do \\
\>\>\>\> calculate the singular pair $(\vs_l,\bdy_l)$ of $
\hat{S}_l(u_{m-1})$ \\
\>\>\>\>\> by iteration (\ref{invit}) \\
\>\>\>\> if $\vs_l < \theta \bnorm{\bdu_{m-1}}$ then \\
\>\>\>\>\> set up the GCD system (\ref{gcdsys}) with $f = u_{m-1}$ \\
\>\>\>\>\>\> (see Section \ref{sec:gcdsys}~) \\
\>\>\>\>\> extract $v_{m}^{(0)}, w_{m}^{(0)}$ from $\bdy_l$ and calculate
$u_{m}^{(0)}$ \\
\>\>\>\>\>\> (see Section \ref{ini_it})~~ \\
\>\>\>\>\>apply the Gauss-Newton iteration (\ref{G-N_it1}) 
from\\ 
\>\>\>\>\>\> $u_{m}^{(0)}, v_{m}^{(0)}, w_{m}^{(0)} $ to obtain $u_{m},
v_{m}, w_{m}$\\
\>\>\>\>\> extract the  residual $\rho = \rho_m$ as in (\ref{gnres_m}) \\
\>\>\>\> end if \\
\>\>\>end do \\
%\>\>\>calculate the (simple) roots of $v_m(x)$\\
\>\>\> adjust the residual tolerance $\varrho$ to be $\max\{\varrho,\;
\phi \rho_j\}$, and  \\
\>\>\>\>\> set $d_m = deg(v_m)$ \\
\>\>end do \\
\>\> set $k=d_1$, $\ell_j = \max \left\{\,t\,\Big|\, d_t \geq k-j+1\,
\right\}, \;\;\; j = 1, 2, \cdots, k.$ \\

\>\> match the roots of $v_m(x)$, $m=1,2,\cdots,s$  \\
\>\>\>\> according to the multiplicities $\ell_j$'s.
\> 
\end{tabbing}
}}}
\end{center} \vspace{-4mm}
\caption{Pseudo-code of Algorithm II}
\label{zgcdroot}
\end{figure}

\subsection{Remarks on the convergence of Algorithm II} \label{lim}

There are two iterative components in Algorithm II.
One of them is the inverse iteration (\ref{invit}).
By Lemma \ref{lem:invit}, the iteration converges for 
all starting vectors $\bdx_0$, unless $\bdx_0$ is orthogonal to the 
intended singular vector $\bdy$. 
The probability of the occurrence of this orthogonality is zero. 
But even if it occurs, roundoff errors in the numerical computation will
quickly destroy the orthogonality during iteration. 
%
%%%%%%%% change back
Therefore, the inverse iteration 
(\ref{invit}) always converges at least slowly. 
The other iterative component is the Gauss-Newton iteration (\ref{G-N_it1})
whose local convergence is ensured in Theorem \ref{gnconv2}. 
Therefore, as long as the rank decision on the Sylvester matrices is
accurate and the error on the initial approximation of the GCD triplet 
is small, Algorithm II will produce correct multiplicity structure 
and a root approximation.

However, due to the nature of the problem, there is no guarantee
that the original multiplicity structure can be identified from an 
inexact polynomial.
When a polynomial is perturbed to a place that has equal distances to 
two or more different pejorative manifolds, it is somewhat
unrealistic to expect any method 
to recover reliably from the perturbation.
Therefore, we have conducted extensive numerical experiments in 
addition to the results exhibited in this paper. 
As reported in our software release note 
\cite{zeng_multroot}, we made a comprehensive test suit of 104 polynomials
based on Jenkins-Traub Testing Principles \cite{jenkins-traub}. 
These polynomials include all the test examples we have seen in the 
literature that have been used by experts to test robustness, 
stability, accuracy and efficiency of root-finders intended for multiple 
roots. 
On all the polynomials with multiple roots in the test suit,
our package {\sc MultRoot} consistently outputs
accurate root/multiplicity results near machine precision. They are
far beyond the ``attainable accuracy'' barrier that other 
algorithms are subject to. 
The report \cite{zeng_multroot} along with the test suit is electronically 
available from the author.
%in the 
%author's homepage\footnote{\tt http://www.neiu.edu/$~$zzeng/multroot.htm}. 

\subsection{Numerical results for Algorithm II }
\label{sec:nrI}

The effectiveness of Algorithm II can be shown by the polynomial
\begin{equation} \label{4rts}
 p(x) = (x-1)^{20}(x-2)^{15}(x-3)^{10}(x-4)^{5} \end{equation}
\begin{wrapfigure}[11]{l}{3.0in}
\vspace{-2.0mm}
\centerline{
\epsfig{figure=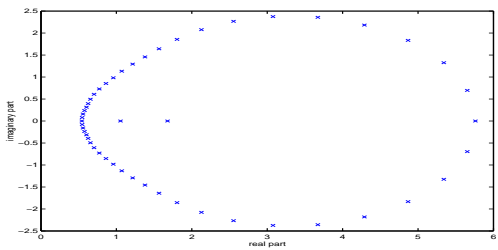, width=2.6in,height=1.27in}} 
\caption[MPSolve results for the polynomial 
(\ref{4rts}) using multiprecision]{\footnotesize
MPSolve results for the polynomial
(\ref{4rts}) using multiprecision}
\label{p1234}
\end{wrapfigure}
generated by Matlab polynomial generator {\sc poly}, with
coefficients rounded up at 16 digits.
Using the default control parameters, Algorithm II code {\sc GcdRoot} 
correctly identifies the multiplicity structure. 
The roots are approximated to an accuracy of 10 digits or more.
With this result as input to Algorithm I code {\sc PejRoot},
we obtained all multiple roots in the end with at least 14 correct 
digits (Table \ref{InII}).

\begin{table}[ht]
\scriptsize
\begin{verbatim}

  Algorithm II (code GcdRoot) result:       | Algorithm I (code PejRoot) result
                                            |   THE BACKWARD ERROR:                  6.16e-016
  The backward error is 6.057721e-010       |   THE ESTIMATED FORWARD ROOT ERROR:    9.46e-014
                                            |
     computed roots       multiplicities    |      computed roots         multiplicities
                                            |
     4.000000000109542           5          |      3.999999999999985               5
     3.000000000176196          10          |      3.000000000000011              10
     2.000000000030904          15          |      1.999999999999997              15
     1.000000000000353          20          |      1.000000000000000              20
\end{verbatim}
\normalsize \vspace{-3mm}
\caption{Roots of $p(x)$ in (\ref{4rts}) computed in two stages}
\label{InII}  \vspace{-4mm}
\end{table}

Polynomials with such high multiplicities are extremely difficult by any 
standard for root-finding. 
The magnitude of its coefficients stretches from 1 to $10^{21}$. 
Remarkably, our algorithms have no difficulty finding all its multiple 
roots. 
To the best of our knowledge, there are no other methods that can calculate 
multiple roots for such polynomials. 
Since the coefficients are inexact, multiprecision root-finders also fail
to calculate the roots with meaningful accuracy. 
Figure \ref{p1234} shows the computed roots by MPSolve 
\cite{bini-mpsolve} using virtually unlimited number of digits in 
machine precision. 
Those results are quite remote from the roots 1, 2, 3, 4.

The Euclidean method has also been used to find GCD in order to identify
the multiplicities \cite{brugnano-trigiante,uhlig-99}.
Uhlig's ~{\sc pzero} \cite{uhlig-99} is a Matlab implementation 
based on the Euclidean method.
The drawback of the Euclidean method is its reliance on recursive long division 
that is numerically unstable (see \S \ref{sec:lsqdiv}). 
%
%Moreover, it is difficult to identify a zero remainder numerically in the 
%process. 
%
%To identify multiplicities, {\sc pzero} must rely on a root 
%cluster matching that is inherently unreliable.
%
Here we compare our code {\sc GcdRoot} with {\sc pzero} on the polynomials
\[ p_k(x) = (x-1)^{4k}(x-2)^{3k}(x-3)^{2k}(x-4)^k
\mbox{ \ \ for $k = 1, 2, \cdots, 8$}. 
\]
When the multiplicities increase, the root accuracy deteriorates with
{\sc pzero}, 
which successfully identifies the multiplicity structure 
for $k=1$ and $k=2$ but fails to do so afterwards.
%
%When $k > 5$, {\sc pzero} no longer recognizes any multiplicities. 
%
%It only outputs root clusters.
%
In comparison, {\sc GcdRoot} consistently attains 
at least 11 digits in root accuracy with increasing multiplicities. 
The multiplicity structures are identified correctly for $k$ up to 7 
and multiplicities up to 28. 
For the current implementation, the limitation of {\sc GcdRoot} on this 
sequence is for $k \leq 7$, whereas the root accuracy will stay the 
same for even larger $k$.
%
%When {\sc GcdRoot} errs at $k=8$, the multiplicities are up 
%to $32$, and the magnitude of coefficients stretches from $1$ to $10^{35}$. 

\begin{table}[ht]
\scriptsize %\tt
\begin{center}
\begin{tabular}{|c|l|lr|lr|lr|} \hline
&&&&&&& \\
{\normalsize $k$} & {\normalsize code} 
& \multicolumn{2}{|c|}{\normalsize $x_1=1$} & 
\multicolumn{2}{|c|}{\normalsize $x_2=2$} & 
\multicolumn{2}{|c|}{\normalsize $x_3=3$} 
%& \multicolumn{2}{|c|}{\normalsize $x_4=4$}
\\ 
&&&&&&& \\ \hline \hline
{\scriptsize $1$} & {\sc pzero}
& 1.00000000001 & (4) & 1.99999999998 &(3) & 3.000000000005 &(2)
%& 4.0000000000004 &(1) 
\\ & {\sc GcdRoot}
& 0.999999999999990 &(4) & 1.99999999999998 &(3) & 3.0000000000005 &(2)
%& 3.9999999999998 &(1) 
\\ \hline
{\scriptsize $2$} & {\sc pzero}
& 1.0000000001 & (8) & 2.000000002 &(6) & 3.000000004 &(4)
%& 4.000000001 &(2) 
\\ & {\sc GcdRoot}
& 0.9999999999998 &(8) & 1.999999999983 &(6) & 2.99999999991 &(4)
%& 3.99999999994 &(2) 
\\ \hline
{\scriptsize $3$} & {\sc pzero}
& 0.9999999897 & ({\bf 13}) & 1.99999990 &({\bf 8}) & 2.9999998 &({\bf 5})
%& 3.99999993 &({\bf 2}) 
\\ & {\sc GcdRoot}
& 0.9999999999997 &(12) & 1.99999999997 &(9) & 2.9999999998 &(6)
%& 3.99999999990 &(3) 
\\ \hline
{\scriptsize $4$} & {\sc pzero}
& 0.9999995 & ({\bf 21}) & 1.999994 &({\bf 6}) & 2.999990 &({\bf 7})
%& 3.999998 &({\bf 2}) 
\\ & {\sc GcdRoot}
& 1.0000000000003 &(16) & 2.00000000002 &(12) & 3.0000000001 &(8)
%& 4.00000000009 &(4) 
\\ \hline
{\scriptsize $5$} & {\sc pzero}
& 1.0000009 & ({\bf 28}) & 2.00001 &({\bf 8}) & 3.00002 &({\bf 6})
%& 4.000004 &(8) 
\\ & {\sc GcdRoot}
& 1.0000000000004 &(20) & 2.00000000003 &(15) & 3.0000000002 &(10)
%& 4.0000000001 &(5) 
\\ \hline
{\scriptsize $6$} & {\sc pzero}
& $----$ & ({\bf 1}) & $----$ & ({\bf 1}) & $----$ &({\bf 1})
%& $----$ &({\bf 1}) 
\\ & {\sc GcdRoot}
& 1.0000000000002 &(24) & 2.00000000001 &(18) & 3.00000000004 &(12)
%& 4.00000000002 &(6) 
\\ \hline
{\scriptsize $7$} & {\sc pzero}
& $----$ & ({\bf 1}) & $----$ & ({\bf 1}) & $----$ &({\bf 1})
%& $----$ &({\bf 1}) 
\\ & {\sc GcdRoot}
& 1.0000000000001 &(28) & 2.00000000001 &(21) & 3.00000000006 &(14)
%& 4.00000000003 &(7) 
%\\ \hline
%{\scriptsize $8$} & {\sc pzero}
%& $----$ & ({\bf 1}) & $----$ & ({\bf 1}) & $----$ &({\bf 1})
%%& $----$ &({\bf 1}) 
%\\ & {\sc GcdRoot}
%& 1.0000000000002 &({\bf 47}) & 2.00000000002 &({\bf 15}) & 
%3.00000000001 &({\bf 11})
%& 4.00000000006 &({\bf 7}) 
\\ \hline
\end{tabular} \end{center}
\normalsize \vspace{-4mm}
\caption[Comparison between {\sc pzero} and {\sc GcdRoot}]{\footnotesize
Partial results on $p_k(x) = (x-1)^{4k}(x-2)^{3k}(x-3)^{2k}(x-4)^k$ and 
comparison between {\sc pzero} and {\sc GcdRoot}.
Numbers in parenthesis are computed multiplicities. Wrong multiplicities
are in boldface.}
\label{table:ulig} \vspace{-5mm}
\end{table}

\section{Numerical results for the combined method} \label{sec:nr}

\subsection{The effect of inexact coefficients}

In application, input data are expected to be inexact. 
The following experiment tests the effect of data error 
on the accuracy as well as robustness of both Algorithm I and II.
For 
\[ p(x) = 
\left(\,x-\frac{\mbox{\scriptsize 10}}{\mbox{\scriptsize 11}} \,\right)^5
\left(\,x-\frac{\mbox{\scriptsize 20}}{\mbox{\scriptsize 11}} \,\right)^5
\left(\,x-\frac{\mbox{\scriptsize 30}}{\mbox{\scriptsize 11}} \,\right)^5
\]
in general form, every coefficient is rounded up to $k$-digit accuracy,
where $k = 10, 9, 8, \cdots$. 

\begin{table}[ht]
\scriptsize
\begin{center}
\begin{tabular}{|c|c|c|l|l|l|c|} \hline
number of & control & &&&&  \\
correct & parameters & code 
& $x_1 = 0.\dot{9}\dot{0}$ 
& $x_2 = 1.\dot{8}\dot{1}$ 
& $x_3 = 2.\dot{7}\dot{2}$ & backward \\
digits & $\varrho$, $\theta $ & & & & & error \\ \hline \hline
$k=10$ & $\varrho = 1e-9$ & {\sc GcdRoot} 
& 0.90909090 & 1.8181818 & 2.7272727 & 1.7e-08 \\
    & $\theta = 1e-7$       & {\sc PejRoot} 
& 0.909090909 & 1.81818181 & 2.7272727 &  2.4e-10 \\ \hline
$k=9$ & $\varrho = 1e-8$ & {\sc GcdRoot} 
& 0.909090 & 1.81818 & 2.72727 & 7.0e-06 \\
    & $\theta=1e-6$    & {\sc PejRoot} 
& 0.9090909 & 1.8181818 & 2.727272 &  2.3e-09 \\ \hline
$k=8$ & $\varrho =1e-7$ & {\sc GcdRoot} 
& 0.90909 & 1.8182 & 2.727 & 1.3e-04 \\
    & $\theta = 1e-5$  & {\sc PejRoot} 
& 0.9090909 & 1.818181 & 2.72727 &  2.3e-08 \\ \hline
$k=7$ & $\varrho=1e-6$ & {\sc GcdRoot} 
& 0.9090  & 1.82   & 2.7 & 1.3e-02 \\
    &  $\theta=1e-4$   & {\sc PejRoot} 
& 0.90909   & 1.81818 & 2.7272  &  2.3e-07 \\ \hline
%$k=6$ & $----$ & {\sc GcdRoot} 
%& $----$    & $----$   & $----$ & $----$    \\
$k=6$ & $----$ & {\sc PejRoot} 
& 0.9090    & 1.8181  & 2.727   &  3.7e-06 \\ \hline
%$k=5$ & $----$ & {\sc GcdRoot} 
%& $----$    & $----$   & $----$ & $----$    \\
$k=5$& $----$ & {\sc PejRoot} 
& 0.909    & 1.818  & 2.72   &  2.4e-05 \\ \hline
%$k=4$ & $----$ & {\sc GcdRoot} 
%& $----$    & $----$   & $----$ & $----$    \\
$k=4$ & $----$ & {\sc PejRoot} 
& 0.90    & 1.81  & 2.7   &  1.9e-04 \\ \hline
%$k=3$ & $----$ & {\sc GcdRoot} 
%& $----$    & $----$   & $----$ & $----$    \\
 $k=3$& $----$ & {\sc PejRoot} 
& 0.9    & 1.8  & 2.8   &  1.8e-03 \\ \hline
\end{tabular}
\end{center}
\normalsize \vspace{-4mm}
\caption{Effect of coefficient error on computed roots}
\label{lowacc} \vspace{-4mm}
\end{table}

For this sequence of problems, Algorithm II code {\sc GcdRoot} correctly
identifies the multiplicity structure if the coefficients have at least 7 
accurate digits. 
If the multiplicities are manually given rather than
computed by {\sc GcdRoot},
Algorithm I code {\sc PejRoot} continues to converge even when data 
accuracy is down to 3 digits. 
For lower data accuracy, the residual tolerance $\varrho$ in {\sc GcdRoot} 
needs to be adjusted accordingly. Table \ref{lowacc} shows the results
of both programs. 

As shown in this test, both methods allow inexact coefficients to 
certain extent. As usual, Algorithm I is more robust than 
Algorithm II, but Algorithm I depends on a structure identifier.

\subsection{The effect of nearby multiple roots}

When two or more multiple roots are nearby, it can be difficult to
identify the correct multiplicity structure. We test the example
\[ p_\eps(x) = (x-1+\eps)^{20}(x-1)^{20}(x+0.5)^5 \]
for decreasing root gap $\eps = 0.1,\; 0.01, \cdots$, making
the root $x_1 = 0.9$, $0.99$, $0.999$, $\cdots$~ along with fixed roots
$x_2=1$ and $x_3=-0.5$.
When root gap decreases, the control parameters may need adjustment.
In this test, we use the default parameters for all cases
except for $\eps = 0.0001$, in which case, the residual growth factor
$\phi = 5$. {\sc GcdRoot} is used to find the initial input for 
{\sc PejRoot}. Computing results are shown for both programs
in Table \ref{rootgap}.

\begin{table}[ht]
\scriptsize
\begin{center}
\begin{tabular}{|l|c|l|l|l|c|r|} \hline
gap        &  & & & & backward & cond. \\
\ \ $\eps$ & code  & $x_1 = 1-\eps $ & 
$x_2 = 1$ & $x_3 = -0.5$ & error & num. \\ \hline \hline
&  {\sc GcdRoot}
& 0.89999999999  & 0.99999999999 & -0.49999999999999 & 9.7e-10 & \\
$10^{-1} $     &         {\sc PejRoot}
& 0.9000000000000 & 0.9999999999999 & -0.50000000000000 &  2.7e-13  &
 .7 \\ \hline
& {\sc GcdRoot}
& 0.98999999  & 0.99999999    & -0.50000000000000 & 3.2e-07 & \\
$10^{-2} $     &         {\sc PejRoot}
& 0.989999999999  & 1.000000000000  & -0.49999999999999 &  1.0e-12 &
6.7 \\ \hline
& {\sc GcdRoot}
& 0.99900     & 1.00000       & -0.49999999999999 & 1.9e-04 & \\
$10^{-3} $     &         {\sc PejRoot}
& 0.99899999999   & 1.00000000000   & -0.500000000000000 &  4.1e-13 &
62.5 \\ \hline
%$\eps= 10^{-4} $ &  {\sc GcdRoot}
%& 0.99994999  & 0.99994999    & -0.5000000000     & 5.7e-08 & \\
%    &        {\sc PejRoot}
%& 0.999949999     & 0.999949999     & -0.500000000       &  2.2e-08 \\ \hline
&  {\sc GcdRoot}
& 0.9997      & 0.99996       & -0.4999999999999  & 1.1e-02 & \\
$10^{-4} $     &        {\sc PejRoot}
& 0.999900000     & 0.999999999     & -0.50000000000000  &  4.0e-12 & 621.7 
\\ \hline
$10^{-5} $ &           {\sc PejRoot} &
 0.999989990     & 1.0000000     & -0.50000000000000  &  4.0e-10 & 5791.8
\\ \hline
\end{tabular}
\end{center}
\normalsize  \vspace{-4mm}
\caption[Effect of decreasing root gap on computed roots]{\footnotesize 
Effect of decreasing root gap on computed roots}
\label{rootgap} \vspace{-3mm}
\end{table}

When the default growth factor stays the same as the default $\phi = 100$
and the gap $\eps \leq 0.0001$, 
{\sc GcdRoot} outputs a multiplicity structure $[40,5]$. 
Namely, {\sc GcdRoot} treats the two nearby 20-fold roots $1$ and $1-\eps$
as a single 40-fold one.
From the computed backward error and the condition number, 
this may not necessarily be incorrect. See Table
\ref{rootgap2}. When backward error becomes $10^{-12}$ and 
condition number is tiny (0.0066), they are 
numerically accurate! In contrast, using the ``correct'' multiplicity
structure $[20,20,5]$, {\sc PejRoot} outputs roots 
with backward error $10^{-10}$ and a large condition 
number $5791.8$ (last line in Table \ref{rootgap}). 

\begin{table}[ht]
\scriptsize
\begin{center}
\begin{tabular}{|l|c|l|l|l|c|r|} \hline
root gap        &  & & & & backward & cond. \\
\ \ \ \ $\eps$ & code  & $x_1 = 1-\eps $ & 
$x_2 = 1$ & $x_3 = -0.5$ & error & num. \ \\ \hline \hline
$\eps=0.0001 $ & {\sc GcdRoot}
& 0.99994999  & 0.99994999    & -0.5000000000     & 5.7e-08 & \\
    &         {\sc PejRoot}
& 0.999949999     & 0.999949999     & -0.500000000       &  2.2e-08  &
0.0066 \\ \hline
$\eps=0.00001 $   & {\sc GcdRoot}
& 0.9999949999 & 0.9999949999  & -0.500000000000   & 1.1e-10 & \\
    &         {\sc PejRoot}
& 0.99999499999   & 0.99999499999   & -0.50000000000     &  4.0e-12 &
0.0066 \\ \hline
\end{tabular}
\end{center}
\normalsize \vspace{-4mm}
\caption[Effect of tiny root gap]{\footnotesize 
If the control parameter is not adjusted, tiny root gap makes computed roots 
identical. 
However, from the backward errors and the condition number, 
they are not necessarily wrong answers.}
\label{rootgap2} \vspace{-3mm}
\end{table}

By adjusting the control parameters, {\sc GcdRoot} can find different
pejorative manifolds that are close to the given polynomial. {\sc PejRoot}
then calculates corresponding pejorative roots. 
The selection of the most suitable solution should be application dependent.

\subsection{A large inexact problem}

\begin{wrapfigure}[16]{r}{1.4in}
\vspace{-3mm}
\scriptsize
\begin{verbatim}
  coefficients of f
      1
     -0.7
     -0.19
      0.177
     -0.7364
     -0.43780
     -0.952494
     -0.2998258
     -0.00322203
     -0.328903811
     -0.4959527435
     -0.9616679762
      0.4410459281
      0.1090273141
      0.6868094008
      0.0391923826
      0.0302248540
      0.6603775863
     -0.1425784968
     -0.3437618593
      0.4357949015
\end{verbatim}
\normalsize
\label{p20c}
\end{wrapfigure}
Implementing the combination of two methods, we have produced a Matlab code 
{\sc MultRoot}. 
We conclude this report by testing this code on our final test problem. 
First of all, twenty complex numbers are randomly generated and
used as roots
\footnotesize
\newline
\(
.5 \pm i,\; 
-1 \pm .2 i,\; 
-.1 \pm i, \; 
-.8 \pm .6i,\;
-.7 \pm .7i,\; 
1.4,\; 
-.4 \pm .9i,\; \)
\newline \(
.9,\; 
-.8 \pm .3i,\;
.3 \pm .8i,\; 
.6 \pm .4i
\)
\normalsize \newline
to generate a polynomial $f$ of degree $20$. 
We then round all coefficients to 10 decimal digits. 
The coefficients are shown in the right. 
We construct multiple roots by squaring $f$ repeatedly. 
Namely, \begin{center}
\( g_k(x) = \left[ \, f(x) \, \right]^{2^k}, \;\;\; k = 1, 2, 3, 4, 5. \)
\end{center}
At $k=5$, $g_5$ has a degree 640 and twenty complex roots of multiplicity 
32. 
Since the machine precision is 16 digits, the polynomials $g_k$
are inexact.  Using the default control parameters, 
our combined program encounters no difficulty in calculating all the roots
as well as finding accurate multiplicities. The worst 
accuracy of the roots is 11-digit. 
Here is the final result.

\vspace{2mm}
\scriptsize
\begin{verbatim}
 THE STRUCTURE PRESERVING CONDITION NUMBER:            0.0780464
 THE BACKWARD ERROR:                         6.38e-012
 THE ESTIMATED FORWARD ROOT ERROR:           9.96e-013
          
 computed roots                   multiplicities  |  computed roots                multiplicities
                                                  |
 0.499999999999399 + 1.000000000006247 i    32    |  1.400000000000303 + 0.000000000000000 i   32
 0.499999999999399 - 1.000000000006247 i    32    | -0.399999999999482 + 0.899999999996264 i   32
-1.000000000003141 + 0.200000000004194 i    32    | -0.399999999999482 - 0.899999999996264 i   32
-1.000000000003140 - 0.200000000004193 i    32    |  0.899999999996995 - 0.000000000000000 i   32
-0.099999999996612 + 1.000000000001018 i    32    | -0.799999999987544 + 0.299999999995441 i   32
-0.099999999996612 - 1.000000000001018 i    32    | -0.799999999987544 - 0.299999999995441 i   32
 0.800000000001492 + 0.600000000001814 i    32    |  0.299999999995789 + 0.799999999976189 i   32
 0.800000000001492 - 0.600000000001815 i    32    |  0.299999999995789 - 0.799999999976189 i   32
-0.699999999997635 + 0.699999999997984 i    32    |  0.599999999989084 + 0.399999999997279 i   32
-0.699999999997635 - 0.699999999997984 i    32    |  0.599999999989084 - 0.399999999997279 i   32
\end{verbatim}
\normalsize

{\bf Acknowledgment. } 
The author wishes to thank the following scholars for their contributions
that improved this paper. 
T. Y. Li, Ross Lippert, Hans Stetter, Joab Winkler and anonymous 
referees made valuable suggestions on the presentation. 
One of the referees pointed out some important previous works in 
\cite{corless-gianni,yun},
Barry Dayton found an error in an early version of the manuscript.
Frank Uhlig provided his insight on the subject in 
e-correspondance along with his software.  
Peter Kravanja also freely shared his code.
D. A. Bini and G. Fiorentino, as well as S. Fortune made their 
root-finders freely available for electronic download.
The author is grateful to the Program Committee of ACM 2003 
International Symposium on Symbolic and Algebraic Computation 
(ISSAC) for their recognition of this work with Distinguished 
Paper Award. 

\bibliographystyle{amsplain}

\begin{thebibliography}{10}

\bibitem{bailey}
{D.~H. Bailey}, {\em A {F}ortran-90 based multiprecision system}, ACM
  Trans. Math. Software, 21 (1995), pp.~379--387.

\bibitem{bini-mpsolve}
{D.~Bini and G.~Fiorentino}, {\em Numerical computation of polynomial roots
  using {MPSolve} -- version 2.0}.
\newblock manuscript, Software and paper available at {\tt
  ftp://ftp.dm.unipi.it/pub/mpsolve/}, 1999.

\bibitem{brugnano-trigiante}
{L.~Brugnanao and D.~Trigiante}, {\em Polynomial roots: the ultimate
  answer?}, Linear Alg. and Its Appl., 225 (1995), pp.~207--219.

\bibitem{brugn}
{L.~Brugnano}, {\em Numerical implementation of a new algorithm for
  polynomials with multiple roots}, J. Difference Eq. and Appl., 1 (1995),
  pp.~187--207.

\bibitem{chin-corless}
{P.~Chin, R.~M. Corless, and G.~F. Corless}, {\em Optimization strategies
  for the approximate {GCD} problem}, Proceedings of 1998 International 
Symposium on Symbolic and Algebraic Computation (ISSAC '98), 
New York, 1998, ACM Press, pp.~228--235.

\bibitem{corless-gianni}
{R.~M. Corless, P.~M. Gianni, B.~M. Trager, and S.~M. Watt}, {\em The
  singular value decomposition for polynomial systems}, 
Proceedings of 1995 nternational
Symposium on Symbolic and Algebraic Computation (ISSAC '95), 
ACM Press, New York, 1995, pp.~195--207.

\bibitem{dedieu-shub}
{J.-P. Dedieu and M.~Shub}, {\em Newton's method for over-determined system
  of equations}, Math. Comp., 69 (1999), pp.~1099--1115.

\bibitem{demmel-ill}
{J.~W. Demmel}, {\em On condition numbers and the distance to the nearest
  ill-posed problem}, Numer. Math., 51 (1987), pp.~251--289.

\bibitem{demmel-kagstrom}
{J.~W. Demmel and B.~Kagstr\"om}, {\em The generalized {S}chur
  decomposition of an arbitrary pencil {$A-\lambda B$}: robust software with
  error bounds and applications. Part {I} \& Part {II}}, ACM Trans. Math.
  Software, 19 (1993), pp.~161--201.

\bibitem{dennis-schnabel}
{J.~E.~Dennis and R.~B. Schnabel}, {\em Numerical Methods for Unconstrained
  Optimization and Nonlinear Equations}, Prentice-Hall Series in Computational
  Mathematics, Prentice-Hall, Englewood Cliffs, New Jersey, 1983.

\bibitem{edelman-elmroth-kagstrom_1}
{A.~Edelman, E.~Elmroth, and B.~Kagstr\"{o}m}, {\em A geometric approach to
  perturbation theory of matrices and and matrix pencils. Part {I}: {V}ersal
  deformations}, SIAM J. Matrix Anal. Appl., 18 (1997), pp.~693--705.

\bibitem{edelman-elmroth-kagstrom_2}
{A.~Edelman, E.~Elmroth, and B.~Kagstr\"{o}m}, {\em A geometric approach
  to perturbation theory of matrices and and matrix pencils. Part {II}: a
  stratification-enhanced staircase algorithm}, SIAM J. Matrix Anal. Appl., 20
  (1999), pp.~667--699.

\bibitem{emiris-galligo-lombardi}
{I.~Z. Emiris, A.~Galligo, and H.~Lombardi}, {\em Certified approximate
  univariate {GCD}s}, J. Pure Appl. Algebra, 117/118 (1997), pp.~229--251.

\bibitem{far-lou-77}
{M.~R. Farmer and G.~Loizou}, {\em An algorithm for the total, or partial,
  factorization of a polynomial}, Math. Proc. Camb. Phil. Soc., 82 (1977),
  pp.~427--437.

\bibitem{far-lou-85}
{M.~R. Farmer and G.~Loizou}, {\em Locating multiple
  zeros interactively}, Comp. Math. Appl., 11 (1985), pp.~595--603.

\bibitem{fortune}
{S.~Fortune}, {\em An iterated eigenvalue algorithm for approximating roots
  of univariate polynomials}, J. Symbolic Comput., 33 (2002), pp.~627--646.

\bibitem{Gautschi-84}
{W.~Gautschi}, {\em Questions of numerical condition related to
  polynomials}, in MAA Studies in Mathematics, Vol. 24, Studies in Numerical
  Analysis, G.~H. Golub, ed., USA, 1984, The Mathematical Association of
  America, pp.~140--177.

\bibitem{hribernig-stetter}
{V.~Hribernig and H.~J. Stetter}, {\em Detection and validation of clusters
  of polynomial zeros}, J. Symb. Comput., 24 (1997), pp.~667--681.

\bibitem{igarashi-ypma}
{M.~Igarashi and T.~Ypma}, {\em Relationships between order and efficiency
  of a class of methods for multiple zeros of polynomials}, J. Comput. Appl.
  Math., 60 (1995), pp.~101--113.

\bibitem{jenkins-traub}
{M.~A. Jenkins and J.~F. Traub}, {\em Principles for testing polynomial
  zerofinding programs}, ACM Trans. Math. Software, 1 (1975), pp.~26--34.

\bibitem{kahan72}
{W.~Kahan}, {\em Conserving confluence curbs ill-condition}.
\newblock Technical Report 6, Computer Science, University of California,
  Berkeley, 1972.

\bibitem{karmarkar-lakshman}
{N.~K. Karmarkar and Y.~N. Lakshman}, {\em On approximate polynomial
  greatest common divisors}, J. Symb. Comput., 26 (1998), pp.~653--666.

\bibitem{kra-van}
{P.~Kravanja and M.~{Van Barel}}, {\em Computing Zeros of Analytic
  Functions, Lecture Notes in Mathematics, 1727}, Springer-Verlag, 2000.

\bibitem{lippert-edelman}
{R.~A. Lippert and A.~Edelman}, {\em The computation and sensitivity of
  double eigenvalues}, in Advances in computational mathematics, Lecture Notes
  in Pure and Appl. Math. 202, New York, 1999, Dekker, pp.~353--393.

\bibitem{miyakoda}
{T.~Miyakoda}, {\em Iterative methods for multiple zeros of a polynomial by
  clustering}, J. Comput. Appl. Math., 28 (1989), pp.~315--326.

\bibitem{pan96}
{V.~Y. Pan}, {\em Numerical computation of a polynomial gcd and
  extensions}.
\newblock Research Report 2996, 
Institut National de Recherche en Informatique
  et en Automatique (INRIA), Sophia-Antipolis, France, 1996.

\bibitem{victorpan97}
{V.~Y. Pan}, {\em Solving polynomial
  equations: some history and recent progress}, SIAM Review, 39 (1997),
  pp.~187--220.

\bibitem{rupprecht}
{D.~Rupprecht}, {\em An algorithm for computing certified approximate {GCD}
  of n univariate polynomials}, J. Pure and Appl. Alg., 139 (1999),
  pp.~255--284.

\bibitem{stetter-cond}
{H.~J. Stetter}, {\em Condition analysis of overdetermined algebraic
  problems}, in Computer Algebra in Scientific Computing--CASC 2000, e.~a.
  V.G.~Ganzha, ed., Springer, 2000, pp.~345--365.

\bibitem{stolan}
{J.~A. Stolan}, {\em An improved \v{S}iljak's algorithm for solving
  polynomial equations converges quadratically to multiple zeros}, J. Comput.
  Appl. Math., 64 (1995), pp.~247--268.

\bibitem{uhlig-99}
{F.~Uhlig}, {\em General polynomial roots and their multiplicities in
  ${O}(n)$ memory and ${O}(n^2)$ time}, Linear and Multilinear Algebra, 46
  (1999), pp.~327--359.

\bibitem{vanhuffel}
{S.~{Van Huffel}}, {\em Iterative algorithms for computing the singular
  subspace of a matrix associated with its smallest singular values}, Linear
  Alg. Appl., 154-156 (1991), pp.~675--709.

\bibitem{wilkinson-63}
{J.~H. Wilkinson}, {\em Rounding Errors in Algebraic Processes},
  Prentice-Hall, Englewood Cliffs, N.J., 1963.

\bibitem{winkler-01}
{J.~R. Winkler}, {\em Condition numbers of a nearly singular simple root of
  a polynomial}, Appl. Numer. Math.,  (2001), pp.~275--285.

\bibitem{ypma}
{T.~J. Ypma}, {\em Finding a multiple zero by transformations and
  {N}ewton-like methods}, SIAM Review, 25 (1983), pp.~365--378.

\bibitem{yun}
{D.~Y.~Y. Yun}, {\em On square-free decomposition algorithms}, in
  Proceedings of 1976 ACM Symposium of Symbolic and Algebraic Computation
  (ISSAC'76), ACM Press, Yorktown Heights, New York, 1976,
  pp.~26--35.

\bibitem{zeng_mult} 
{Z.~Zeng}, \textit{A method computing multiple roots of inexact 
polynomials}, Proceedings of 2003 International Symposium of Symbolic and 
Algebraic Computation (ISSAC '03), 
ACM Press, New York, 2003, pp.~266--272.  

\bibitem{zeng_multroot}
Z.~Zeng, 
{\em Multroot -- a Matlab package computing polynomial roots and
  multiplicities},
\newblock ACM Trans. Math. Software, to appear.

\bibitem{zeng-98}
Z.~Zeng, {\em On ill-conditioned
  eigenvalues, multiple roots of polynomials, and their accurate computation}.
\newblock MSRI Preprint No. 1998-048, (1998). 

\end{thebibliography}

\end{document}